\documentclass[12pt]{report}	

\usepackage{utdiss2}  		

\usepackage{amsmath,amsthm,amsfonts,amscd} 

\usepackage{amssymb}

\usepackage{faktor} 
\usepackage{nicefrac} 
\usepackage{cite} 
\usepackage{enumerate} 

\usepackage{eucal} 	 	
\usepackage{verbatim}      	
\usepackage{makeidx}       	
\usepackage{citesort}         	%
\usepackage{url}		

\author{Mart\'in Mereb}  	

\address{Av. Belgrano 2449 1 G \\ Buenos Aires, Argentina (C1096AAB)}  

\title{On the $E$-polynomials of a family of Character Varieties}

\supervisor
	{Fernando Rodr\'iguez Villegas}

\committeemembers
	[Felipe Voloch]
	[Jeffrey Vaaler]
	[David Ben-Zvi]
	{Gustav Lehrer}

\previousdegrees{Lic.}
\oneandonehalfspacequote

\newtheorem{thm}{Theorem}[section]
\newtheorem{cor}[thm]{Corollary}
\newtheorem{lem}[thm]{Lemma}
\newtheorem{prop}[thm]{Proposition}

\theoremstyle{definition}
\newtheorem{defn}{Definition}[section]

\theoremstyle{remark}
\newtheorem{rem}{Remark}[section]
\newtheorem*{notation}{Notation}
\newtheorem{example}[rem]{Example:}		

\newcommand{\latexe}{{\LaTeX\kern.125em2%
                      \lower.5ex\hbox{$\varepsilon$}}}

\chardef\bslash=`\\	

\makeatletter		
\def\square{\RIfM@\bgroup\else$\bgroup\aftergroup$\fi
  \vcenter{\hrule\hbox{\vrule\@height.6em\kern.6em\vrule}%
                                              \hrule}\egroup}
\makeatother		

\makeindex    

\newcommand{\ip}[2]{\left\langle#1,#2\right\rangle} 
\newcommand{\bra}[1]{\left\langle #1 \right\rangle} 
\newcommand{\pare}[1]{\left( #1 \right)} 
\newcommand{\size}[1]{\left| #1 \right|} 
\newcommand{\set}[1]{\left\{ #1 \right\}} 

\newcommand{\dual}[1]{\wh{ #1 }} 

\newcommand{ \quottt }[2]{{ #1 }/{ #2 }}
\newcommand{ \quott }[2]{\faktor{ #1 }{ #2 } }
\newcommand{ \quot }[2]{\nicefrac{ #1 }{ #2 } }

\newcommand{\mat}[4]{$\begin{pmatrix} #1 & #2 \\ #3 & #4 \end{pmatrix}$} 

\def\noi {\noindent}
\def\Ker {{\rm Ker}}

\def\N {\mathbb{N}}
\def\Z {\mathbb{Z}}
\def\Q {\mathbb{Q}}

\def\C {\mathbb{C}}
\def\F {\mathbb{F}}
\def\S {\mathbb{S}}

\def\k {\mathbb{K}}

\def\l {\lambda}
\def\L {\Lambda}

\def\H {\mathcal{H}}

\def\PP {\mathcal{P}}
\def\Pos{\mathit{P}}

\def\X {\mathfrak{X}}

\def\JO {\Xi}
\def\Ji {\Theta}
\def\Id {Id}
\def\Co {C^0}
\def\delanta {\varphi}
\def\central {\zeta}

\def\MM {\mathcal{M}^g}
\def\MMM {\mathfrak{M}^g}
\def\SSl {\mathfrak{Sl}}

\def\NSol {\widetilde{N}^g}

\def\NSMF {N^g}

\def\d {\delta}
\def\fle {\longrightarrow }

\def\= {\cong}

\def \Spec {\operatorname{Spec}}
\def \gr {\operatorname{gr}}
\def \Aut {\operatorname{Aut}}
\def \Hom {\operatorname{Hom}}
\def \Gl {\operatorname{Gl}}
\def \PGl {\operatorname{PGl}}
\def \Sl {\operatorname{Sl}}

\def \tr {\operatorname{tr}}
\def \Irr {\operatorname{Irr}} 

\def\acts {\curvearrowright}
\def\vec {\overrightarrow}

\def\ol {\overline}
\def\wh {\widehat}

\def\dhat {\wh{d}}

\def\that {\wh{t}}

\def\addsymbol #1: #2#3{$#1$ \> \parbox{5in}{#2 \dotfill \pageref{#3}}\\}

\begin{document}


%
%
%
\commcertpage           

\titlepage              


\begin{dedication}
\index{Dedication@\emph{Dedication}}%
To the $q$.
\end{dedication}
%
%

%
\utabstract
\index{Abstract}%
\indent
We compute the $E$-polynomials of a family of twisted 
character varieties ${\MM}(\Sl_n)$ by proving they have 
polynomial count, and applying a result of N. Katz 
on the counting functions. 

To compute the number of $\F_q$-points of these varieties 
as a function of $q,$ we used a formula of Frobenius.
Our calculations made use of the character tables of 
$\Gl_n(q)$ and $\Sl_n(q),$ previously computed by J. A. 
Green and G. Lehrer, and a result of Hanlon on the 
M\"obius function of a subposet of set-partitions.

The Euler Characteristics of the ${\MM}(\Sl_n)$ are 
calculated then with these polynomial.

%

\tableofcontents   

\listoftables      
\newpage
\chapter*{List of Notation\hfill} \addcontentsline{toc}{chapter}{List of Notation}
\begin{tabbing}
$\zeta_n$~~~~~~~~~~~~~~~~~\=\parbox{5in}{a primitive $n$th-root of unity\dotfill \pageref{e:SMN}}\\
\addsymbol  {[y,x]}: {the commutator of $x$ and $y$, $yxy^{-1}x^{-1}$}{e:SMN}

\addsymbol  {\PP,\PP_n}: {sets of partitions, and partitions of size $n$}{d:Pn}
\addsymbol  {\PP^\JO,\PP_n^\JO}: {sets of $\JO$-partitions, multi-partitions}{n:multip}
\addsymbol  {\L'}: {conjugate multipartition of $\L'$}{n:conjugate}
\addsymbol  {\square,h(\square)}: {a box in a Ferrers diagram, its hook-length}{n:box}
\addsymbol  {\size{\L}}: {size of the multipartition $\L$}{n:multip}
\addsymbol  {H_\l, H_\L,H_\tau}: {Hook polynomials}{e:hookpolyl}
\addsymbol  {\H_\l,\H_\L,\H_\tau}: {Normalized hook polynomials}{n:normhookpoly}
\addsymbol {(\tau)_{d,\l}}: {multiplicity of $(d,\l)$ in $\tau$}{n:typetaudlambda}
\addsymbol {m_{d,\l},m_{d,\l}(\L)}: {multiplicity of $(d,\l)$ in $\L$}{n:multiplicities}
\addsymbol {\tau/(t_d,t_m)}: {quotient type of $\tau$ by $t_d$ an $t_m$}{n:quotype}

\addsymbol  {A \hookrightarrow B}: {injective}{e:NGam}
\addsymbol  {A \twoheadrightarrow B}: {surjective}{n:surj}
\addsymbol  {x \mapsto y}: {$x$ is mapped to $y$}{n:mapsto}
\addsymbol  {A|_{K}}: {restriction of $A$ to $K$}{n:rest}
\addsymbol  {{}^TN}: {transpose of $N$}{e:NGam}

\addsymbol {\vec{\l}}: {a list $(\l_1,\l_2,\ldots,\l_m)$}{e:ZtoS}
\addsymbol {\mu(\;,\;)}: {M\"obius function of a poset}{d:mu}
\addsymbol  {\Pi_n}: {set-partitions of $\set{1,2,\ldots,n}$ }{n:setpart}
\addsymbol  {\Pi_n^\rho,\Pi^\rho}: {set-partitions of fixed by $\rho$}{ex:setp}

\addsymbol  {\bra{\delta}}: {Group generated by $\delta$}{n:deltagen}
\addsymbol  {\dual G}: {dual, or character group, of $G$, $\Hom(G,\C^\times)$}{n:dual}
\addsymbol  {\Irr(G)}: {set of irreducible characters or representatinos of $G$}{n:Irr}


\addsymbol  {Frob_q}: {Frobenius automorphism $x\mapsto x^q$}{n:frob}

\addsymbol  {G \acts X}: {$G$ acts on $X$}{n:newdeg}
\addsymbol  {\faktor{X}{G}}: {set of $G$-orbits of $X$}{r:bijclif}
\addsymbol  {X^\rho}: {$\rho$-fixed points of $X$}{n:xtog}
\addsymbol  {Stab(x)}: {stabilizer of $x$}{t:grnw}
\addsymbol  {Gx}: {$G$-orbit of $x$}{r:tdhat}
\addsymbol  {[x]}: {orbit of $x$}{n:class}
\addsymbol  {\set{\gamma}}: {$Frob_q$-orbit of $\gamma$}{n:froborbit}
\addsymbol  {\deg(\alpha)}: {degree of $\alpha$, the size $\size{\set{\alpha}}$ of its $Frob_q$-orbit}{n:degree}
\addsymbol  {\deg^t(\alpha)}: {newdegree of $\alpha$ relative to $t$, $\deg([\alpha])$ for $[\alpha]\in\faktor{\Gamma}{\bra{\delta^s}}$ }{n:newdeg}
\addsymbol  {(G)_d}: {elements of $G$ whose degree divides $d$, $G^{Frob_q^d}$}{r:weirdgroup}

\addsymbol  {\Sigma_g}: {a genus $g$ Riemann Surface}{c:TQFT}

\addsymbol  {X//G}: {Geometric Invariant Theory quotient of $X$ by $G$}{e:SMN}
\addsymbol  {h^{p,q:j}, h_c^{p,q:j}}: {Hodge numbers}{n:hnum}
\addsymbol  {H(x,y,t;X)}: {Hodge polynomial}{n:hpol}
\addsymbol  {H_c(x,y,t;X)}: {compactly supported Hodge polynomial}{n:hpolc}
\addsymbol  {E(x,y;X)}: {$E$-polynomial, $H_c(x,y,-1;X)$}{n:hpolc}
\addsymbol  {E(q;X)}: {$E$-polynomial, $E(\sqrt{q},\sqrt{q};X)$}{n:Epol}

\addsymbol  {\Gamma}: {Colimit of character groups of $\F_q^\times$}{e:Gamdef}

\addsymbol  {\Ji}: {Frobenius orbits of $\Gamma$, $\faktor{\Gamma}{Frob_q}$}{n:Ji}

\addsymbol {\sigma_\alpha, \tau_\alpha}: {$\dual \F_q$-action by $\alpha$}{n:sigmatau}

\addsymbol {{\MM}(\Sl_n)}: {Twisted Character Variety for $\Sl_n$}{e:SMN}
\addsymbol {{\MM}(\Sl_n)(q)}: {Twisted Character Variety for $\Sl_n(q)$, ${\MM}(\Sl_n(q))$}{e:SMQ}
\addsymbol {\NSMF_n(q)}: {Number of $\F_q$-points of ${\MM}(\Sl_n)$}{n:nsmf}
\addsymbol {\NSol_n(q)}: {Number of solutions of $\displaystyle{ \prod_{i=1}^g[A_i,B_i]=\zeta_n \Id}$}{n:nsol}

\addsymbol  {f(q)|_{q=1}}: {The evaluation of $f$ at $q=1$}{n:eval}
\addsymbol  {\displaystyle \sum_{\substack{i,j\\c(i,j)}}f(i,j)}: {Sum of $f$ over $(i,j)$ satisfying condition $c(i,j)$}{f:mainf}
\end{tabbing} \clearpage

%
%
\chapter{Introduction}
\index{Introduction%
@\emph{Introduction}}%

The aim of this work is to compute the $E$-polynomial of a family of twisted character varieties 
for the group $\Sl_n.$ 
This is accomplished by counting the number $\NSMF_n(q)$\label{n:nsmf} 
of points that these varieties have over finite fields, 
following similar ideas as those from \cite{VilHaus}.

The Character Varieties we study are the affine GIT quotients 
\begin{multline}\label{e:SMN}
 \MM(\Sl_n(\C)):=\{A_1,B_1,\ldots,A_g,B_g\in \Sl_n(\C)\:| \\
 [A_1,B_1]\ldots [A_n,B_n]= \zeta_n I_n\}//\PGl_n(\C) 
\end{multline}
where $[A,B]=ABA^{-1}B^{-1}$, $\zeta_n$ is a primitive $n$-th root of $1$, and the group $\PGl_n(\C)$ acting 
by conjugation.  
They appeared naturally in \cite{VilHaus},(2.2.13) to show that the $\PGl_n$-Character Varieties have an orbifold structure.

The $E$-polynomial of a variety $X$ is 
\begin{equation}\label{e:Epoldef}
E(q ; X ) := H_c ( \sqrt{q},\sqrt{q}, -1; X ) 
\end{equation}
where 
\begin{equation}\label{e:Hpoldef}
H_c( x,y,t ; X):=\sum h_c^{ p,q; j} (X) x^p y^q t^j , 
\end{equation}
with the $h^{p,q;j}$ being the Mixed Hodge numbers of $X,$ 
defined by Deligne in~\cite{Del:2},~\cite{Del:3}.


The strategy to compute the $E$-polynomials of $\MM(\Sl_n)$ is to prove that such varieties have polynomial count, 
which means that there is a polynomial $E_n(q)\in\Z[q]$ such that $\# \MM(\Sl_n(\F_q))=E_n(q)$ for sufficiently 
many prime powers $q$, in the sense described in \ref{s:spread}. 

According to Katz's result from the Appendix of \cite{VilHaus}, the counting polynomials $E_n$ must coincide with 
the $E$-polynomial of the variety.

To compute $\# \MM(\Sl_n(\F_q))$ we regard its $\F_q$-points as: 
\begin{multline}\label{e:SMQ}
 \MM(\Sl_n(\F_q)) = \{A_1,B_1,\ldots,A_g,B_g\in \Sl_n(\F_q)\:| \\
 [A_1,B_1]\ldots [A_n,B_n]= \zeta_n I_n\}//\PGl_n(\F_q) 
\end{multline}
assuming $n|q-1$, a necessary condition for a primitive $n$th root of unity $\zeta_n$ in $\F_q$ to exists.

The problem is then reduced (Corollary \ref{c:sizesorbits}, and Remark \ref{r:numberpoints}) 
to finding an expression for $\NSMF_n(q),$ defined as the number $\NSol_n(q)$ of solutions of the equation
\begin{equation}\label{e:eqQ}
[A_1,B_1]\ldots [A_n,B_n]= \zeta_nI_n ,\; A_i,B_i\in \Sl_n(\F_q) 
\end{equation}
divided by $|\PGl_n(\F_q)|.$

The number of solutions of an equation like \eqref{e:eqQ} over the finite group $G=\Sl_n(\F_q)$ 
can be calculated thanks to the following formula of Frobenius 
\begin{equation}
\label{f:frob}
\sum_{\chi\in \Irr(G)} \left(\frac{|G|}{\chi(1)}\right)^{2g-1} \chi(z) 
\end{equation}
where the sum is taken over all irreducible characters $\chi$ of $G$ (Proposition \ref{p:counting}).

With this end in mind, we need the character table of $\Sl_n(\F_q),$ which was computed 
in~\cite{Leh} and described in~\cite{Green:new} by a technique known as Clifford Theory,
as components of restrictions of irreducible characters of $\Gl_n(\F_q).$ These 
were already calculated in~\cite{Green:old}.

Thanks to these character tables, formula \eqref{f:frob} 
and Katz's Theorem, we prove the following:
\begin{thm} \label{t:main}
The Character Varieties $\MM(\Sl_n)$ have polynomial count and their $E$-polynomials are given by:
\begin{equation} 
E(q; \MM(\Sl_n))=\sum_{\tau,t}\pare{q^{\frac{n^2}{2}}\frac{\H_{\tau'}(q)}{q-1} }^{2(g-1)}t^{2g-1}C_\tau^t.
\end{equation}
\end{thm}
\begin{rem}
The sum is taken over types $\tau$ of size $n$ (defined in \ref{d:type}) and divisors $t$ of $n.$ 
The polynomials $q^{\frac{n^2}{2}}\frac{\H_{\tau'}(q)}{q-1}$ are described in \ref{r:integercoeff}
and the coefficients $C_\tau^t$ are defined in \eqref{e:Ctaut} and computed in \eqref{f:C1toCtaut}.
\end{rem}

Manipulating the terms in this summation we get the following
\begin{cor} \label{c:mainf}
\begin{equation}\label{f:mainf}
E(q; \MM(\Sl_n))=\frac{1}{(q-1)^{2g-2}} \sum_{\substack{|\tau|=n \\ \tau/(t_d,t_m)=\wh{\tau} \\ (\dhat;t_dt_m)=1} }
 \mu(t_d)\frac{O^{2g}(t_dt_m)}{t_dt_m} \pare{ q^{\frac{n^2}{2}}\H_{\tau'}(q) }^{{2g-2}} \Co_{\wh{\tau}}.
\end{equation}
\end{cor}
\begin{rem}
In last formula, $\displaystyle{O^{2g}(t)=\sum_{d|t} \mu\pare{t/d}d^{2g}}$ is the number of elements 
in $(\C^\times)^{2g}$ of order $t,$ the coefficients $\Co_{\tau}$ are $C_\tau/(q-1)$ for the $C_\tau$ 
computed in \cite{VilHaus} 
(see \eqref{f:factor} for a formula), 
$t_d$ and $t_m$ are divisors of $n$ and the sum is taken over a well described set 
of tuples $(\tau,\wh{\tau},t_d,t_m)$ 
(see \ref{f:C1toNSM} for details).
\end{rem}

\begin{rem}
The polynomials found in \ref{f:C1toNSM} 
give the correct number of points $\NSMF_n(q)$ for those $q$ satisfying
\begin{equation}\label{oddity}
q \equiv \left\{
\begin{array}{cc}
1 \mod(n) & \text{when $n$ is odd} \\
1 \mod(2n) & \text{when $n$ is even.}  \end{array}
\right.
\end{equation}
Such stronger condition on $q$ for $n$ even is necessary. For instance, for $n=2$ 
the counting function $\NSMF_2(q)$ is a quasi-polynomial for those $q\equiv 1 \mod(n)$ of period at most $2,$
as we will see in Chapter \ref{ch:example}.
\end{rem}

From these results, and the equidimensionality of the ${\MM}(\Sl_n)$ we get the following:
\begin{cor} \label{c:main}
The $E$-polynomial of $\MM(\Sl_n)$ is palindromic and monic. In particular, the Character Varieties $\MM(\Sl_n)$ are connected. 
\end{cor}

As pointed out in \ref{r:Echar}, one can compute the Euler characteristic of the $\MM(\Sl_n)$ by evaluating their 
$E$-polynomials at $q=1,$ getting
\begin{cor} \label{c:main2}
The Euler Characteristic of the Character Varieties $\MM(\Sl_n)$ is $1$ for $g=1$ and $\mu(n)n^{4g-3}$ for $g\geq 2$.
\end{cor} 



\begin{rem}
The degree of the $E$-polynomial of ${\MM}(\Sl_n)$ is $2(g-1)(n^2-1).$ 
Since the character variety is equidimensional 
$2(g-1)(n^2-1)$ must thus be its dimension. This, together with Corollary \ref{c:main2} 
are consistent with the results of \cite{VilHaus}, more precisely with Corollary 1.1.1 on the Euler 
Characteristic of the $\PGl_n$-Character Varieties and the claims along the proof of Theorem 2.2.12 
on the orbifold structure $${\MM}(\Sl_n) \twoheadrightarrow {\MM}(\PGl_n).$$
\end{rem}


The dissertation is organized as follows: 

In Chapter \ref{ch:prelim} we go over the basics of Combinatorics and Representation Theory that is going to be needed.
In Chapter \ref{ch:charvar} we define the Character Varieties we want to study and
their $E$-polynomials.
In Chapter \ref{ch:example} we elaborate upon the case $n=2$ to show how irreducible characters of $\Gl_2(q)$ split 
and to illustrate the way calculations with formula \eqref{f:frob} carry over.
In Chapter \ref{ch:epoly} we outline the computation of the number of $\F_q$-points of ${\MM}(\Sl_n)$ by applying 
the M\"obius Inversion Formula several many times.
In Chapter \ref{ch:Kitchen} we compute in detail some claims that were left unproven, finishing with the calculation of
$\NSMF_n(q)$ and hence that of $E(q;{\MM}(\Sl_n)).$

\chapter{Preliminaries}\label{ch:prelim}
\index{Preliminaries%
@\emph{Preliminaries}}%

We are going to need some definitions and basic results from Combinatorics and Representation Theory. 
For definitions and notation on partitions and multi-partitions we follow \cite{VinRy}. 
The Frobenius formula for counting solutions of equations on finite groups is explained in \cite{TopGal}.
For the basics on Clifford Theory we follow \cite{Green:new}.

\section{Combinatorics}
\index{Combinatorics%
@\emph{Combinatorics}}%

\subsection{Partitions}
\index{Partitions%
@\emph{Partitions}}%

\begin{notation}\label{d:Pn}
For a partition $\l=(\l_1,\l_2,\ldots,\l_l)$ we have $\l_1\geq\ldots\geq\l_l\geq1$ and write $l(\l)$ for its length $l,$ 
and $|\l|$ for its size $\sum \l_i.$

Let $\PP_n$ note the set of partitions of $n$ and $$\PP:=\bigcup_{n\geq 0}\PP_n.$$ 
We consider the empty partition $\varnothing$ as having size and length $0$, so that $\PP_0$ has one element.
\end{notation}

\begin{rem}
We usually regard $\l_i=0$ whenever $i>l(\l)$ for the ease of notation in formulae. 
\end{rem}

\begin{defn}
For a partition $\l$ we note $\l'$ its conjugate partition $(\l'_1,\ldots,\l'_{l'})$, where $\l'_i$ 
is the number of parts of $\l$ not smaller than $i$.
In particular, $l'=l(\l')=\l_1$ and $l=l(\l)=\l'_1$.
\end{defn}

\begin{rem}
The conjugate partition $\l'$ is usually thought as the one whose Ferrers' diagram (1.3 of \cite{Stan}, for instance)
 is obtained by flipping that of $\l$ across its main diagonal.
\end{rem}

\begin{defn}\label{n:box}
We note $\square\in\l$ a box in the aforementioned diagram, and if $\square$ is in position 
$(i,j)$ then its hook length $h(\square)$ in $\l$ is defined as $$h(\square):= \l_i + \l'_j - i - j +1.$$ 
\end{defn}

\begin{defn}
There is an inner product for partition given by
$${\ip \l \nu}:=\sum_{j\geq 1}\l'_j\nu'_j$$
where 
$\l,\nu\in\PP$, and the number 
$$n(\l):=\frac{1}{2}({\ip \l \l}-|\l|)=\sum_{i\geq 1}(i-1)\l_i.$$
\end{defn}

\begin{example}
For the partition $\l=(7,5,5,3,1)$, we have $\l'=(5,4,4,3,3,1,1)$,
$|\l|=|\l'|=21$ and the Ferrers diagram will be:

\[
\begin{array}{cccccccc}
{}_{\l_i}\diagdown^{\l'_i}& 5 & 4 & 4 & 3 & 3 & 1 & 1 \\ 
7 & \square & \square & \square & \square & \square & \square & \square \\
5 & \square & \square & \square & \square & \square &  &   \\
5 & \square & \blacksquare & \blacksquare & \blacksquare & \blacksquare &  &   \\
3 & \square & \blacksquare & \square &  &  &  &   \\        
1 & \square &  &  &  &  &  & \\
\end{array}
\]
the black squares corresponds to the hook of $\square_{3,2}$ the box in position $(3,2)$ and its hook length is:
$$h\left(\square_{3,2}\right)= 5+4 -3-2 +1=5.$$

\end{example}

\subsection{Multi-Partitions}
\index{Multi-Partitions
@\emph{Multi-Partitions}}%

We will consider collections $\JO$ of finite sets $X\in\JO$. In Section \ref{s:charGl} they will be $Frob_q$-orbits of characters.
\begin{defn}\label{n:multip}
For such a collection $\JO$, we define an $\JO$-partition (or multi-partition) as a function $\L:\JO\to\PP$. 

The size of a multi-partition $\L$ being defined as 
$$|\L|:=\sum_{X\in\JO}|X||\L(X)|.$$
Let $\PP^\JO_n$ note the set of multi-partitions of size $n$ and 
$$\displaystyle{\PP^\JO:=\bigcup_{n\geq 0}\PP^\JO_n}.$$
\end{defn}

\begin{defn}\label{n:conjugate}
For $\L \in \PP^\JO_n$ we define its conjugate $\L' \in\PP^\JO_n$ by $$\L'(X):=\L(X)'\in \PP.$$
\end{defn}


\begin{defn}
Similarly we extend the notion of number for multi-partitions as 
$$n(\L):=\sum_{X \in \JO} |X|n(\L(X)).$$ 
\end{defn}

\begin{notation}
We will usually regard $\L\in\PP^\JO$ as a function from $\bigcup_{X\in\JO}X$ to $\PP$ constant on the sets $X.$ 
Thus, for $x\in X,$ $\L(x):=\L(X).$
\end{notation}

\begin{defn}
The support of $\L$ is the union of all those $X\in\JO$ mapped to a nonzero partition $\L(X)\ne \varnothing.$
\end{defn}

\begin{defn}\label{n:multiplicities}
Given an integer $d \geq 1$ and a partition $\l \in \PP$ we define the multiplicity of $(d,\l)$ in $\L$ as 
$$m_{d,\l}=m_{d,\l}(\L):=\#\set{ X \in\JO\:|\: \size{X} = d,\L(X)=\l },$$ 
when $\l\neq 0$, and set $m_{d,0}=0$ as a convention. 
\end{defn}

\begin{defn}\label{d:type}
Let $\tau(\L)$ be the collection of multiplicities $\pare{m_{d,\l}}_{d\geq 1,\l\in\PP}$ and call it the type of $\L.$ 
Its size is given by $\size{\tau}:=\size{\L}=\sum m_{d,\l}d\size{\l}.$
\end{defn}

\begin{notation}
We write $(\tau)_{d,\l }$ or $m_{d,\l}(\tau)$ for the multiplicity $m_{d,\l}(\L)$ of any multi-partition $\L$ of type $\tau.$
\end{notation}

\begin{defn}\label{n:typetaudlambda}
The support of $\tau$ is the set of pairs $(d,\l)$ with a nonzero associated multiplicity $(\tau)_{d,\l }.$
\end{defn}

\begin{defn}
We write $\square \in \L$ for a box in the Ferrers diagram of one of the partitions $\L(X)$ and 
define its hook length $h(\square)$ 
in $\L$ as $|X|$ times its hook length in $\L(X)$.
\end{defn}

\subsection{Hook Polynomials}
\index{Hook Polynomials%
@\emph{Hook Polynomials}}%

\begin{defn}
Given a partition $\l$ we associate to it its hook polynomial 
\begin{equation}\label{e:hookpolyl}
H_\l(q):=\prod_{\square\in\l}(q^{h(\square)}-1). 
\end{equation}

In a similar fashion, for a multi-partition $\L$ we define:
\begin{equation}\label{e:hookpolyL}
H_\L(q):=\prod_{\square\in\L}(q^{h(\square)}-1)=\prod_{X\in\JO}H_{\L(X)}(q^{|X|}).
\end{equation}
\end{defn}

\begin{notation}
Since the $H_\L(q)$  only depends on the type $\tau$ of $\L$ we will write $H_\tau$ for $H_\L.$
\end{notation}

Most formulas will look more transparent if we introduce the normalized version of the hook polynomial.
\begin{defn}\label{n:normhookpoly}
For $\l\in\PP$ and $\L\in\PP_n(\Gamma)$ any map of type $\tau=\{m_{d,\l}\}$, the corresponding normalized hook polynomials are:
\begin{equation}\label{e:hookl}
\H_\l(q):=q^{-\frac{1}{2} \ip{\l}{\l}}\prod_{ \square \in \l }(1-q^{h(\square)}) \in \Z[q^{\frac{1}{2}},q^{-\frac{1}{2}}]
\end{equation}
and
\begin{equation}\label{e:hookL}
\H_\L(q):=\prod_{\{\gamma\}}\H_{\L(\gamma)}(q^{deg(\gamma)})=\prod_{d,\l}\H_\l(q^{d})^{m_{d,\l}}.
\end{equation}
Since $\H_\L$ depends only on the type $\tau$ of $\L,$ we also write $\H_\tau(q) = \H_\L(q).$
\end{defn}

\begin{rem} By the symmetry of the factors of \eqref{e:hookl} we have the identity
\begin{equation}\label{e:palin}
\H_\l(q^{-1})=(-1)^{|\l|}\H_{\l'}(q). 
\end{equation}
\end{rem}

\begin{defn}\label{d:fortypes}
For a type $\tau$ the size, number, length and hook polynomial of $\tau$ are then defined as those of $\L$ 
for any $\L$ of type $\tau$.
In a similar way, we can also define the conjugate $\tau'$ as that of $\L'.$ 
\end{defn}

\subsection{M\"obius Inversion formula:}\label{s:MIV}
\index{M\"obius Inversion formula
@\emph{M\"obius Inversion formula}}%

Let $(\Pos,\geq)$ be a finite poset, i.e. a finite partially ordered set.


\begin{defn} 
To every function $f:\Pos \to \C$ we assign 
\[\wh f :\Pos\to\C\]
\[\wh f (a) := \sum_{b\geq a}f(b) \]
the accumulated sums of $f$ with respect to $\Pos.$
\end{defn}

\begin{rem}\label{n:mapsto}
We see that $f \mapsto \wh f $ is a linear operator in $\Hom_{\C}(\C^\Pos,\C^\Pos)$ whose matrix $M_\Sigma$ is given by
$$(M_\Sigma)_{a,b}=\begin{cases} 1 \mbox{ if } b \geq a\\ 0 \mbox{ otherwise.}\end{cases} $$
the adjacency matrix of (the digraph induced by) $\Pos.$

Since $\Pos$ is a poset, this matrix will be lower triangular with respect to some total ordering refining that 
of $\Pos$, and will have only $1$'s in its main diagonal.

Therefore its inverse $\mu\in\Hom_{\C}(\C^\Pos,\C^\Pos)$ has also only $1$'s in the diagonal and integer entries, for writing $M_{\Sigma}$ as $Id+A$, with $A$ strictly lower triangular (hence nilpotent), the alternating sum $$Id - A +A^2 - A^3 \ldots$$ for $(Id+A)^{-1}$ results finite.
\end{rem}

\begin{defn}\label{d:mu}
The function $\mu(a,b)$ given by the entries $(\mu)_{a,b}$ of the matrix $\mu$ 
is called the M\"obius Function for $\Pos$.
\end{defn}

\begin{rem}
By definition of $\mu$ we have the M\"obius Inversion Formula:
\begin{equation}\label{f:MoIn}
f(a):= \sum_{b\geq a}\mu(a,b){\wh f}(b).
\end{equation}
\end{rem}

%

\begin{example}\label{ex:div}
\emph{(Divisors)} let $\Pos$ be the set of positive divisors of a fixed $n$ with the ordering given by reversed 
divisibility, more precisely $a\geq b$ if and only if $a|b$.

The accumulated sums for $f$ are $\displaystyle{ \wh{f}(m):=\sum_{d|m}f(d)}$ and the inversion formula for this 
case is the well known
\begin{equation}\label{f:mIdiv}
f(m):=\sum_{d|m}\wh{f}(d)\mu(m/d) 
\end{equation}
where 
\begin{equation*}
\mu(m):=
\begin{cases} (-1)^{\#\{\text{primes dividing $m$}\}} & \text{for square-free $m$,}
\\
0 &\text{otherwise.}
\end{cases}
\end{equation*}
\end{example}

In other words $\mu(a,b)= \mu(a/b)$ if $b$ divides $a$ and $0$ if not.
\begin{example}\label{n:setpart}
\emph{(Set-partitions)} Let us take $\Pos$ as the collection $\Pi_n$ of set-partitions of the set 
$\{1,2,\ldots, n\}$ for a fixed $n$ and the order being given by reversed
refinement. 

For instance, in $\Pi_4$ we have
$$ 1234 \geq 12|34 \geq 1|2|34 \geq 1|2|3|4.$$

Then the M\"obius function is:
$$\mu(\nu,\pi)=(-1)^{s-r}\prod_{i=1}^s (i-1)!^{r_i}$$
where $s,r$ are the numbers of parts of $\nu$ and $\pi$ (respectively), $r_i$ is the number of parts 
of $\pi$ containing exactly $i$ parts of $\nu$. 
\end{example}

\begin{example}\label{ex:setp}
\emph{(Fixed set-partitions)} Given a poset $\Pi_n$ as in last example, and a $\rho\in\S_n$ we take the subposet 
$\Pi_n^\rho$ of those set-partitions fixed by $\rho$. In general, its 
M\"obius function is not easy to compute.
In \cite{Han}
Hanlon computed $\mu(\wh 0,\wh 1)$ where $\wh 0$ and $\wh 1$ stands for the bottom and the top element in such a poset, 
namely those partitions with $n$ and $1$ parts, respectively. 
His results reads:
\begin{equation}\label{f:Han}
\mu(\wh 0,\wh 1)= 
\begin{cases} \mu(d)(-d)^{m_d-1}(m_d-1)! & \text{if $\rho$ is a product of $m_d$ $d$-cycles,}
\\
0 &\text{otherwise.}
\end{cases}
\end{equation}
\end{example}

\begin{rem}
Thhe last example played its role in the computation of the $E$-polynomials of the twisted $\Gl_n$-character 
varieties in \cite{VilHaus}.
We are going to use it for the $\Sl_n$-case as well. 
\end{rem}

\begin{rem}\label{r:torus}
Example \ref{ex:div} is also going to be important for many computations, 
as we will see in Chapter \ref{ch:epoly} and Chapter \ref{ch:Kitchen}.

As an application of it, let us compute the number $O^N(n)$ of elements of a order in the torus $(\C^\times)^N$.
Writing $\wh O^N(m)$ for the number of elements in the $m$-torsion subgroup, the following two observations leads to the desired formula:
\begin{itemize}
\item $\displaystyle{\wh O^N(m)= \sum_{d|m} O^N(d)}.$
\item $\displaystyle{\wh O^N(m)= m^N}.$
\end{itemize}
Therefore 
\begin{equation}\label{f:torus}
O^N(n)=\sum_{d|n}\mu\pare{\frac{n}{d}}d^N. 
\end{equation}
\end{rem}


\section{Representation Theory}
\index{Representation Theory
@\emph{Representation Theory}}%


In this section we will list the facts from Representation Theory we will need.
We assume $G$ is a finite group, and our representations will all have finite dimension.
For proofs and details we refer the reader to \cite{Serre}, \cite{FulHar}.

A representation of $G$ is a finite dimensional $\C$-vector space $V$ together with a group homomorphism
$$\rho:G\to \Aut_{\C}(V).$$ 

Its character $\chi=\chi_\rho$ is the class function obtained from composition with the trace 
$$\chi_{\rho}(g):=\tr(\rho(g)).$$

Since the character $\chi_\rho$ of a representation $\rho$ determines it uniquely up to isomorphism, 
a representation is usually identified with its character.

A representation is said to be irreducible if it does not have nontrivial invariant subspaces.
In such case, its character is also called irreducible. 
\begin{notation}\label{n:Irr}
For a finite group $G$ we note $\Irr(G)$ the set of its irreducible characters. Whenever it is clear from the context, 
it will also represents the set of (isomorphism classes of) irreducible representations of $G.$
\end{notation}

In an irreducible representation $(V,\rho)$, Schur's Lemma says that all the endomorphisms $\Hom_\rho(V,V)$ are 
the scalar multiplications $\central \Id: x\mapsto \central x$ with $\central\in\C.$

There is a natural inner product $\ip{\;\: }{\: }$ for characters given by
$$\ip{\chi}{\chi'}:= \frac{1}{|G|}\sum_{g\in G}\chi(g) \ol{\chi'(g)}.$$ 

We also recall the Orthogonality relations:

For $\chi,\chi'\in\Irr(G)$ then 
\begin{equation}\label{e:orto1}
\ip{\chi}{\chi'}=
\begin{cases} 1 \mbox{ if } \chi = \chi' \\ 0 \mbox{ otherwise.}
\end{cases} 
\end{equation}

For $g,h\in G$
\begin{equation}\label{e:orto2}
\sum_{\chi \in \Irr(G)} {\chi(g)}\ol{\chi'(h)}=
\begin{cases} \size{C_g(G)} \mbox{ if $g$ and $h$ are conjugate}  \\ \quad \;0 \quad\;\; \mbox{ otherwise}
\end{cases} 
\end{equation}
where $C_g(G)$ is the centralizer of $g$ in $G$ (2.2 in \cite{FulHar}, for instance).

\subsection{Character Groups}
\index{Character Groups
@\emph{Character Groups}}%

\begin{defn}\label{n:dual}
For $G$ a finite group, we define its character group (or dual) as $\dual G := \Hom(G,\C^\times).$ 
\end{defn}
\begin{rem}
It has a structure of abelian group given by pointwise multiplication. Its identity element is the trivial (or principal) 
character $1_G: g\mapsto 1.$
\end{rem}
\begin{rem}
For $G$ abelian $\dual G = \Irr(G).$ 
\end{rem}
\begin{rem}
When $G$ is finite and cyclic, $\dual G$ is easily seen to be cyclic of same order, then
$G \= \dual G$ in a non canonical way, and the same holds for any finite abelian group $G$.
\end{rem}

\begin{rem}
The canonical map $G\fle \dual{\dual{G}}$ is an isomorphism.
\end{rem}

\begin{rem}\label{r:subquotdual}
Since $\C^\times$ is divisible, any short exact sequence of abelian groups
$$0\fle K \stackrel{i}{\fle} G \stackrel{\pi}{\fle} Q \fle 0 $$
induces its dual
$$0\fle \dual Q \stackrel{\pi^*}{\fle} \dual G \stackrel{i^*}{\fle} \dual K \fle 0 $$
leading us to the conclusion
\begin{itemize}
 \item The dual of a quotient is a subgroup of the dual:
 namely, the dual of $Q = \quott{G}{K}$ is the subgroup of those $\phi\in\dual G$ that vanish over $K$.
 \item The dual of a subgroup is a quotient of the dual:
 namely, the dual of $K\subseteq G$ is the quotient of $\dual G$ by the dual $\dual Q.$
\end{itemize}
\end{rem}

\begin{rem}\label{r:group}
For $G$ cyclic there is only one subgroup and only one quotient of order $d$ 
for every $d$ divisor of $|G|.$
Given $H\subseteq G$ and $g\in H$, $g\in H$ if and only if $|g|$ divides $H.$
\end{rem}

\begin{rem}\label{r:fullcharsum}
For $G$ abelian,
$$
\sum_{\phi \in \dual G} \phi(g) =\begin{cases} |G|\mbox{ if } g=1 \\ \;\: 0\;\: \mbox { otherwise.}\end{cases}
$$
from \eqref{e:orto2}, for instance.
\end{rem}

\begin{rem}\label{r:partialcharsum}
Given a subgroup $\dual Q$ of $\dual G,$ by Remarks \ref{r:fullcharsum} and \ref{r:subquotdual} we have
$$
\sum_{\phi \in \dual Q} \phi(g) =\begin{cases} |Q|\mbox{ if } g\in K \\ \;\: 0\;\: \mbox { otherwise.}\end{cases}
$$ 
\end{rem}

\begin{rem}\label{r:group2}
For cyclic $G$, by Remarks \ref{r:group} and \ref{r:partialcharsum}, 
the computation of a character sum over a subgroup of $\dual G$ is reduced to a question of divisibility of integers. 
This is going to be important in the proof of Lemma \ref{l:nonzero}.
\end{rem}

\subsection{Counting solutions of equations in finite groups}
\index{Counting solutions of equations in finite groups
@\emph{Counting solutions of equations in finite groups}}%

Let $G$ be a finite group and $A$ some abelian group like $\C$ or $\C^{n\times n}$. 
\begin{defn}
For a function $f:G\to A$ we note
$$\int_G f(x)dx := \frac{1}{|G|}\sum_{x\in G} f(x)$$ the integral of $f$ with respect to the Haar Measure of $G$.
\end{defn}

\begin{prop}\label{p:counting}
Given $z\in G$, the number of $2g$-tuples $(x_1,y_1,\ldots,x_g,y_g)$ satisfying 
$[y_1,x_1]\ldots [y_g,x_g]z=1$
is:
\begin{equation}\label{f:counting}
\#\{ [y_1,x_1]\ldots [y_g,x_g]z=1 \} = \sum_{\chi\in\Irr(G)} \chi(z) \left(\frac{|G|}{\chi(1)}\right)^{2g-1}.  
\end{equation}
\end{prop}
\proof

Given $\rho:G\to \Aut(V)$ an irreducible representation of $G$ and $\chi$ its character, let us consider
for a fixed $x\in G$ the average 
\begin{equation}\label{e:average}
\int_G \rho(yxy^{-1})dy. 
\end{equation}
Since it commutes with the $G$-action, it must be, by Schur's lemma, a scalar map $\central \Id$. 

By taking traces we get:
\begin{eqnarray}
\int_G \rho(yxy^{-1})dy & = & \central \Id \nonumber\\
\int_G \tr(\rho(yxy^{-1}))dy & = &\central \tr(\Id) \nonumber\\
\int_G \chi(yxy^{-1})dy & = &\central \chi(1) \nonumber\\
\int_G \chi(x)dy & = &\central \chi(1) \nonumber\\
\frac{\chi(x)}{\chi(1)}& = &\central  \label{e:count1}
\end{eqnarray}
thus, the average \eqref{e:average} becomes
\begin{equation}
\int_G \rho(yxy^{-1}) dy= \frac{\chi(x)}{\chi(1)}\Id. \label{e:count2}
\end{equation}

Multiplying by $\rho(x^{-1})$ from the right we get:
$$\int_G \rho(yxy^{-1}x^{-1}) dy = \frac{\chi(x)}{\chi(1)}\rho(x)^{-1}$$
summing over all $x\in G$ and dividing by $|G|$ again we end up with
\begin{equation}\label{e:count3}
\iint_{G^2} \rho([y,x])dxdy= \int_G \frac{\chi(x)}{\chi(1)}\rho(x)^{-1}dx. 
\end{equation}
Since the left hand side of this equation is invariant under $G$-conjugation, it is also a scalar transformation.
Taking traces again we conclude
\begin{eqnarray}
\iint_{G^2} \rho(yxy^{-1}x^{-1})dxdy & = & \int_G \frac{\chi(x)\chi(x^{-1})}{\chi(1)^2}dx \Id \nonumber \\ 
& = & \int_G \frac{\chi(x)\ol{\chi(x)}}{\chi(1)^2}dx \Id \nonumber \\ 
& = & \frac{1}{\chi(1)^2} \Id.\label{e:count4}
\end{eqnarray}

Raising \eqref{e:count4} to the $g$th power and multiplying by $\rho(z)$ from the right
we will have
\begin{eqnarray}
\left(\iint_{G^2} \rho(yxy^{-1}x^{-1})dxdy\right)^{g}\rho(z) & = & \frac{1}{\chi(1)^{2g}} \rho(z) \nonumber\\ 
\underbrace{\int_{G} \ldots \int_{G}}_{2g \text{ times}} \rho([y_1,x_1]\ldots [y_g,x_g]z)dx_1dy_1\ldots dx_gdy_g  
& = & \frac{1}{\chi(1)^{2g}} \rho(z).\label{e:count5}
\end{eqnarray}
Taking trace at both sides yet one more time we have:
\begin{equation}\label{e:count6}
\int_{G} \ldots \int_{G} \chi([y_1,x_1]\ldots [y_g,x_g]z)dx_1dy_1\ldots dx_gdy_g   
=  \frac{\chi(z)}{\chi(1)^{2g}} . 
\end{equation}

Multiplying this by $\chi(1)$ and summing over all $\chi\in \Irr(G)$ we see that the sum 
\begin{equation}\label{e:righthandside}
\sum_{\chi\in\Irr(G)} \frac{\chi(z)}{\chi(1)^{2g-1}} 
\end{equation} 
is equal to
\begin{equation}\label{e:count7}
\int_{G} \ldots \int_{G} \sum_{\chi\in\Irr(G)}\chi(1)\chi([y_1,x_1]\ldots [y_g,x_g]z)dx_1dy_1\ldots dx_gdy_g 
\end{equation}
and thanks to the orthogonality relations 
from \eqref{e:orto2} only those terms with $[y_1,x_1]\ldots [y_g,x_g]z=1$ survive in \eqref{e:count7}, and we get
\begin{eqnarray}\label{f:countfinale}
\frac{1}{|G|^{2g-1}}\#\{ [y_1,x_1]\ldots [y_g,x_g]z=1 \} &=& \sum_{\chi\in\Irr(G)} \frac{\chi(z)}{\chi(1)^{2g-1}} 
\end{eqnarray} 
from which the proposition follows immediately.
\endproof

\begin{rem}
In Chapter \ref{ch:epoly} we will compute, with the aid of formula \eqref{f:counting}, the 
number of $\F_q$-points in the twisted $\Sl_n$-character varieties, 
mimicking the ideas of \cite{VilHaus} for the $\Gl_n$-case.
\end{rem}

As a particular case of Proposition \ref{p:counting}, for $z=1$ (the identity element of $G$) 
we recover the well known formula (see~\cite{FdQ} for instance).
\begin{cor}\label{c:TQFT}
Let $\Sigma_g$ be a genus $g$ compact Riemann surface, $\pi_1$ its fundamental group and $G$ a finite group. 
Then
\begin{equation}
\frac{1}{ |G|} \#\Hom(\pi_1(\Sigma_g),G)  =  \sum_{\chi\in\Irr(G)} \left(\frac{|G|}{\chi(1)}\right)^{2g-2}. 
\end{equation}
\end{cor}

\subsection{Clifford Theory}\label{s:cliff}
\index{Clifford Theory%
@\emph{Clifford Theory}}%

To apply formula \ref{p:counting}, we will need the character table of $\Sl_n(q),$ computed by Lehrer in \cite{Leh}.
In order to deal with these characters we will also make use of the results of \cite{Green:new}.

The standard setting for Clifford Theory is the following: 

$G$ is a finite group and $H$ is a normal subgroup of $G$ such that the quotient $\quott{G}{H}$ is abelian.

There is a natural action of $G$ on $\Irr(H)$ 
given by the conjugation 
$$(g\theta)(h):=\theta(g^{-1}hg).$$ 
Since the characters of $H$ are class functions, $H$ acts trivially on $\Irr(H),$ therefore this $G$-action 
induces one of $\quott{G}{H}$. 

We note $[\theta]$ the class of $\theta\in\Irr(H)$ under this action.\label{n:class}

There is also a natural action of $\dual{\pare{\quott{G}{H}}}:=\Hom\pare{\quott{G}{H},\C^\times},$ 
the character group of $\quott{G}{H}$, on $\Irr(G)$ given by multiplication 
$$(\psi\chi)(g):=\psi([g])\chi(g)$$
where $\psi \in\dual{\pare{\quott{G}{H}}}$ and $[g]\in{\quott{G}{H}}$ is the class of $g\in G$.

We note $[\chi]$ the orbit of $\chi\in\Irr(G)$ under this action.

\begin{thm} (1 from \cite{Green:new})\label{t:grnw}
Let $\chi,\chi' \in \Irr(G).$ Then:
\begin{itemize}
\item $\chi_H = \chi'_H$ if and only if $\chi'\in  [\chi]$.
\item $ |Stab (\chi)|= \ip{\chi_H}{\chi_H }$.
\item $\chi_H$ is irreducible if and only if $|Stab (\chi)|=1$.
\item if $\theta \in \Irr(H)$ is a constituent of ${\chi_H}$, then there is a positive integer $e_H(\chi)$ such that
$$ \chi_H = e_H(\chi)\sum_{\theta'\in [\theta]}\theta'$$ and $$e_H(\chi)^2 \size{[\chi]}\size{[\theta]} = |G/H|.$$
\item $e_H(\chi) =  1$ for $\quottt{G}{H}$ cyclic.
\end{itemize}
\end{thm}

\begin{rem}\label{r:bijclif}
When the quotient $\quottt{G}{H}$ is cyclic, theorem \ref{t:grnw} establishes a canonical bijection between both sets of orbits 
$$\quott{\Irr(G)}{\wh{\pare{\quottt{G}{H}}}}\qquad \text{        and         }\qquad \quott{\Irr(H)}{\pare{\quottt{G}{H}}} $$ 
by which an orbit $[\chi]$ is mapped to the $\quott{G}{H}$-orbit $[\theta]$ defined as 
the set of irreducible components of $\displaystyle{\chi_H=\sum_{\theta'\in [\theta]} \theta'}$, and 
the product of the sizes of both orbits is $\size{\quott{G}{H}}.$
\end{rem}

\subsection{Characters of finite General Linear Groups}\label{s:charGl}
\index{Characters of finite General Linear Groups
@\emph{Characters of finite General Linear Groups}}%

Thanks to Clifford Theory and motivated by the formulas for the $\PGl_n$-character varieties from \cite{VilHaus}, 
we need a description of the irreducible characters of the finite general linear groups $\Gl_n(q)$. 

They were computed by Green in \cite{Green:old}. 
We will follow Macdonald's approach form \cite{MacD}, with some notation taken from \cite{VinRy}.

Let us fix a prime power $ q = p^m.$ 
We note $ \F_{q^r} $ the field of $ q^r $ elements inside a fixed algebraic closure $ \ol{\F}_q $ of $ \F_q $. 

\begin{notation}\label{n:frob}
Let $ Frob_q : x \mapsto x^q $ be the Frobenius automorphism. 
\end{notation}

\begin{defn}\label{n:xtog}
Clearly $\F_{q^r}$ is $ \ol{\F_q}^{Frob_q^r} $ the field of elements fixed by $Frob_q^r.$ 
\end{defn}

\begin{notation}\label{n:surj}
Whenever $n|m$ we have the norm map $ N^m_n:\F_{q^m}\twoheadrightarrow\F_{q^n} $ which is known to be surjective.
\end{notation}

\begin{defn}\label{d:Gam}
Let $ \Gamma_n:=\Hom(\F_{q^n}^\times,\C^\times) $ the character group of the cyclic group 
$\F_{q^n}^\times.$ By Remark \ref{r:subquotdual}, the transpose of the norm gives us an 
injective map 
\begin{equation}\label{e:NGam}
{}^TN^m_n:\Gamma_n \hookrightarrow \Gamma_m
\end{equation}
when $ n|m $, defining a direct system whose colimit we note 
\begin{equation}\label{e:Gamdef}
\Gamma := \varinjlim\Gamma_n .
\end{equation}
\end{defn}

\begin{rem}\label{r:identify}
The Frobenius automorphism induces by transposition, a map on $ \Gamma $ given by $ \gamma \mapsto \gamma^q $ and 
just as before we can identify $ \Gamma_n $ with $\Gamma^{Frob_q^n},$ the elements of $ \Gamma $ fixed by $ Frob_q^n $. 
\end{rem}

\begin{defn}\label{n:degree}
The degree of $\gamma \in \Gamma$ is $\deg(\gamma)\in \N$ the size of the Frobenius orbit of $\gamma.$ 
\end{defn}

\begin{rem}
Since we are identifying $\Gamma_n$ with elements of the colimit $\Gamma$, given some $\gamma \in \Gamma$ and an element 
$\xi\in\F_q^\times$ there are in general more than one way to interpret the evaluation $\gamma(\xi),$ for if $\gamma$ belongs to 
$\Gamma_n =\Hom(\F_{q^n}^\times,\C^\times)$ there is another representative ${}^TN_n^{nk}(\gamma)$ for $\gamma$ in every $\Gamma_{nk},$ and 
the evaluation of the new representative at $\xi$ gives
${}^TN_n^{nk}(\gamma)(\xi) = \gamma(N_n^{nk}(\xi))=\gamma(\xi^k)= \gamma(\xi)^k.$ 

To avoid confusion we will say in which character group $\Gamma_n$ the character is supposed to be considered.
\end{rem}

\begin{rem}\label{r:evalua}
For instance, given $\gamma\in\Gamma_d=\Hom(\F_{q^d}^\times,\C^\times)$ the product 
$${}^TN_1^d(\gamma):=\gamma \gamma^q \gamma^{q^2}\ldots \gamma^{q^{d-1}}$$
is $Frob_q$-stable, hence there is a representative of it in $\Gamma_1=\Hom(\F_q^\times,\C^\times).$ Let us call such representative 
$\gamma_1\in\Hom(\F_{q}^\times,\C^\times).$

Evaluation at any $\xi\in\F_q^\times$ of ${}^TN_1^d(\gamma)$ in $\Gamma_d$ will be 
$$\gamma(\xi)\gamma(\xi^q) \ldots \gamma(\xi^{q^{d-1}})=\gamma(\xi)^d$$
since $\xi^q=\xi.$

On the other hand, evaluation in $\Gamma_1$ (i.e.: computing $\gamma_1(\xi)$) amounts to evaluating at $\xi$ the restriction 
to $\F_q^\times$ of any of the $Frob_q$-conjugates of $\gamma,$ since all such conjugates agree on $\F_q^\times$ 
and the restriction $\gamma|_{\F_q^\times}\in\Gamma_1$ \label{n:rest}
will become ${}^TN_1^d(\gamma)$ under the inclusion $\Gamma_1 \hookrightarrow \Gamma_d.$
\end{rem}

\begin{defn}\label{n:Ji}
Consider $\Ji:=\quott{\Gamma}{Frob_q}$ the collection of $Frob_q$-orbits of $\Gamma$.
\end{defn}

\begin{rem}\label{n:froborbit}
Every $Frob_q$-orbit $\set{\gamma}$ for $\gamma\in\Gamma$ is finite.
\end{rem}

\begin{rem}
By last remark, the set $ \PP^\Ji $ of $\Ji$-partitions $\L$ is defined, and so are 

\begin{itemize}
\item the size $\displaystyle \size{\L} = \sum_{\gamma \in \Gamma}|\L (\gamma)|,$

\item the conjugate $\L'(\gamma)=\L(\gamma)'\in \PP,$

\item the multiplicities $$m_{d,\l}=m_{d,\l}(\L)=\#\set{\set{\gamma}| \deg(\gamma)=d,\L(\gamma)=\l},$$

\item the type $\tau(\L)= \pare{m_{d,\l} }_{d\geq 1,\l\in\PP},$

\item the hook polynomials $H_\L,$ $H_\tau$ and

\item the set $\PP^\Ji_m$  of multi-partitions size $m$, for $m\geq 0.$
\end{itemize} 
\end{rem}

\begin{thm} \label{t:bijMac} 
((6.8) from \cite{MacD} page 286) 
There is a bijective 
correspondence between $\PP_n^\Ji$ and the irreducible characters of $Gl_n(q)$ by which the character $\chi_\L$ 
corresponding to $\L$ has degree 
\begin{equation}\label{e:degree}
\chi_\L(1)=\frac{\prod_{i=1}^n(q^i-1)}{ q^{-n(\L')} H_{\L}(q)}
\end{equation}
and the value at the central element $\central \Id$ ($\central \in \F_q^*$) is given by 
\begin{equation}\label{e:central}
\chi_\L(\central \Id)=\Delta_\L(\central)\chi_\Lambda(1)
\end{equation}
where 
\begin{equation}\label{e:delta}
\Delta_\L=\prod_{\gamma \in \Gamma} \gamma^{|\L(\gamma)|}\in \Gamma_1 .
\end{equation}
\end{thm}

\begin{notation}
We write $\Delta_\chi:=\Delta_\L $ for $\chi=\chi_\L$, the character $\chi\in\Irr(G)$
associated to the function $\L\in\PP_n^\Ji.$ 
\end{notation}

\begin{rem}\label{r:deltafactorcentral}
Since $\Delta_\L$ is $Frob_q$-stable, the evaluation of \eqref{e:central} should be done as explained in 
Remark \ref{r:evalua}, meaning
\begin{eqnarray}\label{e:deltafactorcentral}
\Delta_\L(\central)=\prod_{\gamma \in \Gamma} \gamma^{|\L(\gamma)|}(\central) \\ \nonumber
 = \prod_{\set{\gamma} } \gamma^{|\L(\gamma)|}(\central)
\end{eqnarray}
where last product is taken over the $Frob_q$-orbits $\set{\gamma}$ and each evaluation $\gamma(\central)$ is done in
$\Gamma_{\deg(\gamma)}.$
\end{rem}

\begin{rem}\label{r:principal}
The trivial character $1_G$ corresponds to the $\Ji$-partition $\L\in\PP^{\Ji}_n$ sending $1\in\Gamma_1$ to $(1,\ldots,1)$ of size $n$ 
(and hence mapping all other characters to $\varnothing$).
\end{rem}

\begin{rem}\label{r:1diml}
More generally, from \eqref{e:degree} follows that all $1$-dimensional representations of $\Gl_n(q)$ 
are in one-to-one correspondence with the multi-partitions 
$\L\in\PP^{\Ji}_n$ of type $\tau_p$ with $(\tau_p)_{1,(1,\ldots,1)}=1.$ 
\end{rem}

\begin{rem}\label{r:integercoeff}
Since the degree of a character $\chi_\Lambda$ depends only on the type 
$\tau=\tau(\Lambda)$ we write $\chi_\tau(1):=\chi_\Lambda(1)$ and from \ref{e:degree} we get
\begin{equation}\label{e:quot}
\frac{|Gl_n(q)|}{\chi_\tau(1)}= (-1)^nq^{\frac{n^2}{2}}\H_{\tau'}(q)
= q^{\binom{n}{2}-n(\L')}H_{\L}(q)
\end{equation}
which actually lies in $\Z[q].$
\end{rem}
%
%
%

\begin{defn}\label{d:circproduct}
The $\circ$-product of characters $\chi_1$ from $\Gl_{n_1}(q)$ and $\chi_2$ from $\Gl_{n_2}(q)$ is defined
by considering the parabolic subgroup $P$ of upper triangular block matrices of the form:
\begin{center}
{\mat{A_1}{\ast}{0}{A_2}} 
\end{center}
with $A_i\in\Gl_{n_i}(q)$ and its character $\chi_1 \ast \chi_2$ constructed from
$
\chi_1\otimes\chi_2 
$ 
by composition with the natural projection
$$P \twoheadrightarrow \Gl_{n_1}(q)\times\Gl_{n_2}(q).$$
Then take $$\chi_1\circ\chi_2 := (\chi_1 \ast \chi_2)^G=Ind_P^G(\chi_1 \ast \chi_2)$$ the induced character.
\end{defn}

\begin{rem}
In \cite{Green:old}, the character $\chi_\L\in\Irr(\Gl_n(q))$ is constructed 
for $\L\in \PP_n^{\Ji}$ as
$$
\chi_\L = J^{\bra{\gamma_1}}(\L(\gamma_1))\circ \dots \circ J^{\bra{\gamma_m}}(\L(\gamma_m))
$$
where $J^{\bra{\gamma}}(\l)$ is certain simple character of $\Gl_{\deg(\gamma)|\l|}(q).$ 
\end{rem}
\begin{rem}
Since the $\circ$-product is easily seen to be associative and commutative, the character 
$J^{\bra{\gamma_1}}(\L(\gamma_1))\circ \dots \circ J^{\bra{\gamma_m}}(\L(\gamma_m))$ 
is well defined.
\end{rem}

\subsection{The actions $\sigma_\alpha$ and $\tau_\alpha$}\label{n:sigmatau}\label{ch:sigmatau}
\index{The actions $\sigma_\alpha$ and $\tau_\alpha$%
@\emph{The actions $\sigma_\alpha$ and $\tau_\alpha$}}%


The determinant 
identifies the quotient $\quott{\Gl_n(q)}{\Sl_n(q)}$ with the cyclic group $\F_q^\times$.

The $\dual{\F}_q^\times$-action on $\Irr(\Gl_n(q))$ we talked about in Section \ref{s:cliff} corresponds to the
$\sigma_\alpha$ operation defined for $\alpha\in\dual{\F}_q^\times$ in \cite{Green:new} as 
$$\sigma_\alpha(\chi)(x):=\alpha(\det(x))\chi(x).$$


There is another action defined in \cite{Leh}. Namely the $\tau_\alpha$ action defined by
$$
\tau_\alpha(\chi) := J^{\bra{\alpha\gamma_1}}(\L(\gamma_1))\circ \dots \circ J^{\bra{\alpha\gamma_m}}(\L(\gamma_m))
$$ 
where the product $\alpha\gamma_i$ should be computed in
$\Gamma_{d_i}$ (following Remark \ref{r:evalua}).

\begin{rem}
The operation $\chi\mapsto \tau_\alpha(\chi)$ from \cite{Leh} 
in the language of $\Ji$-partitions is given by
$$\L \quad \longmapsto \quad \L \alpha^{-1}$$
where $$\L \alpha^{-1}(\gamma):=\L(\alpha^{-1}\gamma).$$
\end{rem}

\begin{rem}\label{r:action}
The proposition of page 133 in \cite{Green:new} asserts that both operations $\tau_\alpha$ and $\sigma_\alpha$ coincide. 
%
This will help us to regard $\Irr(\Gl_n(q))$ as a set of functions, and $\Irr(\Sl_n(q))$ as a set of orbits 
of those functions under the recently described action. 
\end{rem}

Theorem \ref{t:grnw} and last Remark together imply the following
\begin{thm}[4 in \cite{Leh}, 2 in \cite{Green:new}]
Let $\chi, \chi' \in \Irr(\Gl_n(q))$, then their restrictions 
$\chi_{\Sl_n(q)}$ and $\chi'_{\Sl_n(q)}$ agree if and only if there exists $\alpha\in\dual \F_q^\times$ such 
that $\tau_\alpha(\chi)=\chi'.$ 
\end{thm}

\chapter{The $E$-polynomial of a Character Variety}\label{ch:charvar}
\index{The $E$-polynomial of a Character Variety
@\emph{The $E$-polynomial of a Character Variety}}%


\section{Character Varieties}
\index{Character Varieties
@\emph{Character Varieties}}%

In this Section we define the family of Character Varieties we want to study, and list some of their properties.


\begin{defn}\label{d:SMn}
Let $n$ and $g$ be a positive integers, $\k$ a field possessing a primitive $n$th root of unity we call $\zeta_n,$ and
$\Id\in\Sl_n(\k)$ the identity matrix.
The twisted $\Sl_n$-character variety ${\MM}(\Sl_n(\k))$ is given by all the classes of $2g$-tuples $(A_1,B_1,\ldots,A_g, B_g)$ 
in $\Sl_n(\k)^{2g}$ satisfying the the following matricial equation
\begin{equation}\label{e:SMndef1}
 [A_1,B_1]\ldots[A_g,B_g]=\zeta_n Id
\end{equation}
modulo $\PGl_n(\k)$-conjugation.
In other words, it is the geometric quotient 
\begin{equation}\label{e:SMndef2}
{\MM}(\Sl_n(\k)):= \{ [A_1,B_1]\ldots[A_g,B_g]=\zeta_n Id\} // \PGl_n(\k)
\end{equation}
where $A_i,B_i \in \Sl_n(\k).$
\end{defn}
\begin{rem}
They appeared in Theorem 2.2.12 of \cite{VilHaus}, to prove that the $\PGl_n$-character variety ${\MM}(\PGl_n)$ 
is an orbifold.
\end{rem}

\begin{rem}
For a finite field $\F_q$ to have a primitive $n$th root of unity it is necessary and sufficient that $n$ divides $q-1.$
\end{rem}

\begin{notation}
When $\k=\F_q$ is the finite field with $q$ elements we shall write ${\MM}(\Sl_n)(q)$ for ${\MM}(\Sl_n(\F_q)).$ 
\end{notation}

\begin{rem}
The twisted $\Sl_n$-character varieties are non-singular and equidimensional (see page 570 of \cite{VilHaus}).
\end{rem}

\begin{prop}\label{p:ssimple}
Let $(A_1,B_1,\ldots,A_g, B_g) \in \Sl_n(\k)^{2g}$ satisfying \eqref{e:SMndef1}, and $Z\in \k^{n\times n}$ a matrix 
commuting with all the $A_i$'s and $B_i$'s. Then $Z$ is a scalar matrix (i.e.: a scalar multiple of the identity $Id$).
\end{prop}
\proof 
The only nonzero subspace $K\subseteq \k^{n}$ invariant by the $A_i$'s and $B_i$'s is $\k^{n}$, since restriction to such $K$
will lead to the identity
\begin{equation}\label{e:Subs}
 [A_1|_K,B_1|_K]\ldots[A_g|_K,B_g|_K]=\zeta_n Id_K
\end{equation}
which, by taking determinants, gives $1=\zeta_n^{\dim_{\k}K},$ implying $\dim_\k K = n$ and then $K=\k^{n}.$

Extending $\k$ if necessary, we can assume $Z$ has an eigenvalue $\l\in\k$. Since the eigenspace $K=\Ker(Z-\l Id)$ is going to be invariant
by multiplication by the $A_i$'s and $B_i$'s, we conclude that $K=\k^{n},$ or what is the same, $Z=\l Id.$
\endproof

\begin{rem}
The stabilizer of a $2g$-tuple $(A_1,B_1,\ldots,A_g, B_g) \in \Sl_n(\k)^{2g}$ satisfying \eqref{e:SMndef1} 
under the $\Gl_n(\k)$-action of simultaneous conjugation is, according to proposition \eqref{p:ssimple}, 
the center of $\Gl_n(\k).$  
\end{rem}

Therefore we have the following two corollaries
\begin{cor}\label{c:freeaction}
The induced $\PGl_n(\k)$-action results free.
\end{cor}

\begin{cor}\label{c:sizesorbits}
When $\k=\F_q$, the number of $\k$-points of ${\MM}(\Sl_n)(q)$ will be the number of solutions of \eqref{e:SMndef1} 
in $\Sl_n(q)$ divided by $|\PGl_n(q)|.$
\end{cor}


\section{Mixed Hodge Structures}
\index{Mixed Hodge Structures@\emph{Mixed Hodge Structures}}%

In \cite{Del:2} and \cite{Del:3} Deligne proved the existence of the following structure.
Let X be a complex algebraic variety. 
For each $j$ the cohomology group $H^j(X,\Q)$ and its complexification $H^j(X,\C)$ have filtrations
\begin{equation}\label{e:filtW}
0= W_{-1}\subseteq W_0 \subseteq \ldots \subseteq W_{2j}=H^j(X,\Q) 
\end{equation}
and
\begin{equation}\label{e:filtF}
H^j(X,\C) = F^{0} \supseteq F^{1} \supseteq \ldots \supseteq F^{m}\supseteq F^{m+1}=0 
\end{equation}
called weight and Hodge filtration, respectively, with the property that 
the filtration induced by $F$ on $\gr_l W := W_l /W_{l-1}$ (the complexification of the graded pieces of the
weight filtration) equips every graded piece with a pure Hodge structure of weight $l.$ 

\begin{rem}\label{r:propMHS}
Such a structure is called the Mixed Hodge Structure, 
or MHS, of $X$. It behaves nicely with respect 
to the most common operations on cohomology. 
For instance, it is preserved by maps 
$$f^* : H^*(Y,\Q) \to H^*(X,\Q)$$ 
induced from algebraic maps $f : X \to Y$,
by maps induced from field automorphisms $\sigma \in \Aut( \C/\Q )$, 
by K\"unneth isomorphism 
$$H^*(X\times Y,\Q)\=  H^*(X,\Q)\otimes H^*(Y,\Q)$$ 
and by cup products.

For smooth $X$ we get trivial $W _{j-1} H^j (X)$
.
\end{rem}
\begin{rem}\label{r:poinca}
There is also a Mixed Hodge Structure for cohomology with compact support $H_c^*(X,\Q)$, and the interplay with
$H^*(X,\Q)$ given by Poincar\'e
duality respects both MHS.
\end{rem}

\begin{defn}\label{d:hnumdef}
Define:
\begin{equation}\label{n:hnum}
 h^{ p,q; j} (X) := \dim_\C \gr_p F \gr_{p+q} W H^j(X,\C) , 
\end{equation}
and 
\begin{equation}\label{e:hnum}
h_c^{ p,q; j} (X) := \dim_\C \gr_p F \gr_{p+q} W H_c^j(X,\C)  
\end{equation}
the Mixed Hodge numbers and the compactly supported Mixed Hodge numbers, respectively. 
\end{defn}
These numbers can be encoded in the following polynomials:
\begin{defn}\label{d:hpoldef}
The Mixed Hodge Polynomial
\begin{equation}\label{n:hpol}
 H(x,y,t;X ):=\sum h^{ p,q; j} (X) x^p y^q t^j ,
\end{equation}
the compactly supported Mixed Hodge Polynomial
\begin{equation}\label{n:hpolc}
H_c( x,y,t ; X):=\sum h_c^{ p,q; j} (X) x^p y^q t^j , 
\end{equation}
and finally let
\begin{equation}\label{n:Epol}
E(x, y ; X ) := H_c ( x, y, -1; X ) 
\end{equation}
be the $E$-polynomial of $X$.
\end{defn}

\begin{rem}\label{r:propertiesEpoly}
From the definition of $E(x, y ; X )$ , \eqref{e:hnum},\eqref{e:filtW} and \eqref{e:filtF} we have: 
\begin{itemize}
\item $E(1,1;X)$ is the Euler Characteristic of $X$
\item The total degree of $E(x,y;X)$ is twice the dimension of $X.$
\item The coefficient of $x^{dim(X)}y^{dim(X)}$ in $E(x,y:X)$ is the number of connected components of $X.$
\end{itemize}
\end{rem}

%
%

\section{Katz's Theorem}\label{s:spread}
\index{Katz's Theorem@\emph{Katz's Theorem}}%

\begin{defn}\label{d:spreading}
For a complex algebraic variety $X$, a finitely generated $\Z$-algebra $R$ and a fixed embedding $$\phi:R\hookrightarrow\C$$
we say that a separated scheme $\X/R$ is a spreading out of $X$ if its extension of scalars $\X_\phi$ is isomorphic to $X$.
\end{defn}

\begin{rem}\label{r:spreading}
A spreading out is an object that links a complex algebraic variety $X$ with its finite field counterparts.
\end{rem}

\begin{defn}\label{d:polcount}
Let us now assume $X$ has a spreading out $\X$ such that for every ring homomorphism $\psi:R\to\F_q$, the number of points of
$\X_\psi(\F_q)$ is given by $P_X(q)$ for some fixed $P_X(t)\in\Z[t],$ independent of $\psi$.
We say that this $X$ has polynomial count and that $P_X$ is the counting polynomial.
\end{defn}

The following result holds:
\begin{thm}\label{t:Katz}
[(2.1.8) in~\cite{VilHaus} proved by Katz in Appendix] Let $X$ be a variety over $\C$. Assume $X$ has 
polynomial count with polynomial $P_X(t)\in\Z[t],$ 
then the $E$-polynomial of $X$ is given by $E( x, y ; X) = P_X(xy).$
\end{thm}

\begin{notation}\label{n:Epoly}
In such case, we write 
$$E(q; X ):=E(\sqrt{q},\sqrt{q}; X ).$$
\end{notation}

\begin{rem}\label{r:Epoly}
This polynomial $E(q; X )\in\Z[q]$ agrees with the counting polynomial $P_X(q)$ 
and encodes all the information that $E(x, y; X )$ does. 
\end{rem}

\begin{rem}\label{r:Echar}
For $X$ a variety with polynomial count given by $P_X(q),$ by Remark \ref{r:propertiesEpoly} we have:
\begin{itemize}
\item The Euler characteristic of $X$ can be computed as $P_X(1).$
\item The dimension $dim(X)$ is the degree of $P_X(q).$
\item The principal coefficient of $P_X(q)$ is the number of connected components of $X.$
\end{itemize}
\end{rem}




\section{Spreading out of ${\MM}(\Sl_n)$}\label{s:spreadfinitefield}
\index{Spreading out of ${\MM}(\Sl_n)$@\emph{Spreading out of ${\MM}(\Sl_n)$}}%


Let us fix $\zeta_{2n} \in\C$ and $\zeta_n=\zeta_{2n}^2$ primitive $2n$th and $n$th-roots of $1,$ respectively.

Let $R$ be the finitely generated $\Z$-algebra $Z[\zeta_{2n},\frac{1}{n}]\subseteq \C.$

Consider $\wh{\MMM}(\Sl_n)$ the closed affine subscheme of $\SSl_n^{2g}$ over $R$ defined 
by the matricial equation \eqref{e:SMndef1}.

\begin{defn}\label{d:spreadoutcharvar}
Define ${\MMM}(\Sl_n)$ as the categorical quotient:
\begin{equation}\label{e:sesh}
{\MMM}(\Sl_n):=\Spec(R[\wh{\MMM}(\Sl_n)]^{\PGl_n(R)}).
\end{equation}
\end{defn}

\begin{rem}\label{r:spreadout}
By flatness of $R\hookrightarrow\C$ and Lemma 2 from \cite{Sesh} we have
\begin{equation}\label{e:sesh2}
R[\wh{\MMM}(\Sl_n)]^{\PGl_n(R)} \otimes_R \C = \C[\wh{\MMM}(\Sl_n)]^{\PGl_n(\C)}.
\end{equation}
Then ${\MMM}(\Sl_n)$ is a spreading out of ${\MM}(\Sl_n).$
\end{rem}

\begin{rem}\label{r:spoutfinite}
By Lemma 3.2 from the appendix of \cite{Kac}, the $\wh{\MMM}(\Sl_n)$-fibers of an $\F_q$-point of ${\MMM}(\Sl_n)$ will be nonempty 
$\PGl_n(\F_q)$-orbits, hence
\begin{equation}\label{e:kac}
{\MMM}(\Sl_n)(\F_q)=\wh{\MMM}(\Sl_n(\F_q))/\PGl_n(\F_q)
\end{equation}
and by Corollary \ref{c:freeaction} the action will be free.
\end{rem}

\begin{rem}\label{r:numberpoints}
Last Remark, together with Corollary \ref{c:sizesorbits} imply that the counting functions for our character varieties 
${\MM}(\Sl_n)$ are 
\begin{equation}\label{e:nsmfspread}
\NSMF_n(q):=\frac{1}{|\PGl_n(\F_q)|}\NSol_n(q)
\end{equation}
where
\begin{equation}\label{e:nsol}
\NSol_n(q)=\sum_{\theta\in \Irr(\Sl_n(\F_q))} \pare{ \frac{|\Sl_n(\F_q)|}{\theta(1)} }^{2g-1}\theta(\zeta_n\Id). 
\end{equation}
is the number of solutions of \eqref{e:SMndef1} (thanks to Proposition \ref{p:counting}).
\end{rem}

\begin{rem}\label{r:goodqs}
Note that since $R\subseteq\C$ is defined as the finite $\Z$-algebra $Z[\zeta_{2n},\frac{1}{n}],$
every ring homomorphism to a field $F$
$$\phi: R \to F$$
must necessarily send $n$ to a unit. 

Since 
$$
(1-\zeta_{2n})(1-\zeta_{2n}^2)\ldots(1-\zeta_{2n}^{2n-1})=2n
$$
when $n$ is even, $2n$ must also be mapped to a unit, and hence $\zeta_{2n}$ is mapped to an element of order $2n.$
In particular, a primitive $n$th root of unity $\zeta_n$ is mapped to an element of order $n.$

When $n$ is odd, $2$ could be in the kernel of $\phi$, making $F$ a characteristic $2$ field (in which case no element has even order).
A similar reasoning with
$$
(1-\zeta_{n})(1-\zeta_{n}^2)\ldots(1-\zeta_{n}^{n-1})=n
$$
proves that $\zeta_n$ is mapped to an order $n$ element.

Hence, in both cases we have an element of order $n$ in the target field.

When $F=\F_q$ is a finite field with $q$ elements, $F$ has an element of order $n$ if and only if $n$ divides $q-1.$
When $n$ is even, $F$ must also have an order $2n$ element, making $q-1$ divisible by $2n.$

Summarizing, under this setting, our counting functions only make sense for those $q$ satisfying \ref{oddity}.
\end{rem}

\begin{rem}\label{r:quasigoodqs}
If we took $R=\Z[\zeta_n, \frac{1}{n}]$ instead, most of the above follows, except that the only condition 
for $q$ is to be congruent to $1$ modulo $n.$ This is all what we need for ${\MM}(\Sl_n(\F_q))$ to make sense.
The counting function for this spreading out will not in general be a polynomial (as we will see in Chapter \ref{ch:example}), and 
that is why we take $R$ with a primitive $2n$th root of $1.$
\end{rem}


\chapter{The Case $n=2$}\label{ch:example}
\index{The Case $n=2$%
@\emph{The Case $n=2$}}%


In this Chapter we are going to list the characters of $\Gl_2(q)$ and see how their restrictions to $\Sl_2(q)$ 
split for odd $q,$ giving all the irreducible characters of such subgroup as we saw in \ref{s:cliff}.
We then use these characters to compute the number of points of ${\MM}(\Sl_2)(q)$ for $q$ odd, and hence their $E$-polynomials 
(by Theorem \ref{t:Katz}). We finally get some geometric and topological information thanks to Remark \ref{r:Echar}.

\section{Character Table of $\Gl_2(q)$}
\index{Character Table of $\Gl_2(q)$
@\emph{Character Table of $\Gl_2(q)$}}%

As shown in \cite{Fran} and \cite{FulHar}, 
the character table of $\Gl_2(q)$ is:
\noi
\begin{table}[h] \label{t:one}
\begin{center}
\caption{Character Table of $\Gl_2(q)$.}
\vskip 10pt
\begin{footnotesize}
\begin{tabular}{|c|cccc|}
\hline
 classes & \mat a 0 0 a & \mat a 1 0 a & \mat a 0 0 b & \mat x 0 0 {x^q} \\
 \hline
\# of classes & $q - 1$ & $q - 1 $&$ (q-1)(q-2)/2 $&$ q(q-1)/2 $\\
 class size & $1$ & $q^2 - 1$ & $q(q+1) $&$ q(q-1) $ \\
\hline
$R_{T}^G(\alpha,\beta)$ & $(q+1)\alpha(a)\beta(a)$ & $\alpha(a)\beta(a)$ & $\alpha(a)\beta(b)+ \alpha(b)\beta(a)$ &$ 0 $\\
$-R_{T^s}^G(\omega)$ & $(q-1)\omega(a)$ & $-\omega(a)$ & $0$  & $-(\omega(x)+\omega(x^q))$ \\
$\sigma_\alpha (1_G)$ & $\alpha(a^2)$ & $\alpha(a^2)$ & $\alpha(ab)$ & $\alpha(x^{q+1})$\\
$\sigma_\alpha (St_G)$ & $q\alpha(a^2)$ & $0$ & $\alpha(ab)$ & $-\alpha(x^{q+1})$ \\
\hline
\end{tabular}\\[10pt]
\end{footnotesize}
\vskip -20pt
\end{center}
\end{table}
\index{commands!environments!table}%

\noi
where $a,b\in\F_q^\times$, $a\neq b$, $x\in\F_{q^2}\setminus\F_q$, $\alpha,\beta\in \wh{\F}_q^\times$, $\alpha\neq\beta$,
$\omega\in\wh{\F}_{q^2}^\times$, $\omega\neq\omega^q$, $St_G$ is the Steinberg representation and
$1_G$ the trivial one.

The class \mat x 0 0 {x^q} corresponds to those matrices in $\Gl_2(q)$ having eigenvalues $x$ and $x^q$, 
so that they are diagonalizable in $\F_{q^2}$ but not in $\F_q$.

The rows correspond to the different types $\tau_1$, $\tau_2$, $\tau_3$ and $\tau_4$ of size $2,$
all of them having a unique nonzero multiplicity
$(\tau_1)_{1,(1)}=2,$ $(\tau_2)_{2,(1)}=1,$ 
$(\tau_3)_{1,(2)}=1$ and $(\tau_4)_{1,(1,1)}=1$ respectively.

\begin{itemize}
\item Type $\tau_1$ has two characters $\alpha \neq \beta$ 
of degree 1 mapped to the partition $(1)\in\PP_1$ and the notation $R_{T}^G(\alpha,\beta)$ corresponds to 
Deligne-Lusztig \cite{DelLus}
. We must point out that $R_{T}^G(\alpha,\beta) \= R_{T}^G(\beta,\alpha).$

\item Analogously, type $\tau_2$ (one degree $2$ character $\omega\in\Gamma_2$
mapped to $(1)\in\PP_1$) corresponds to $-R_{T^s}^G(\omega)$ where $T^s$ is the non-split torus (see \cite{DelLus}). 
In this case we have $-R_{T^s}^G(\omega) \= -R_{T^s}^G(\omega^q).$

\item For every $\alpha\in\Gamma_1=\dual \F_q^\times$, the composition with the determinant 
$\det:\Gl_n(q)\to\F_q^\times$ 
gives a new character we noted $\sigma_\alpha(1_G)$ (defined in Subsection \ref{ch:sigmatau}). Type $\tau_3$ 
is the collection of all these.

\item Type $\tau_4$ consists on all the $\sigma_\alpha(St_G)$ found by tensoring one character of 
$\tau_3$ with the Steinberg representation.
\end{itemize}

\section{Character Table of $\Sl_2(q)$}
\index{Character Table of $\Sl_2(q)$
@\emph{Character Table of $\Sl_2(q)$}}%

To compute the character table of $\Sl_2(q)$ we will apply the results from Subsections \ref{s:cliff} and \ref{ch:sigmatau}.
Let us now assume $q$ odd. 

The action of $\dual \F_q^\times$ on $\Irr(\Gl_2(q))$ identifies its orbits with representations of $\Sl_2(q)$ 
(by restriction) each of which is a sum of representations in an orbit of the $\F_q^\times$-action given by 
conjugation (under the identification $\F_q^\times\= \Gl_2(q)/\Sl_2(q)$)
. And the product of the sizes of corresponding orbits must be $q-1$.

\begin{rem}\label{r:irredn3}
The characters of type $\tau_3$ have all dimension $1$ and hence, their restrictions to $\Sl_2(q)$ have to be irreducible.
In fact, all of them restrict to the trivial representation of $\Sl_2(q).$
\end{rem}

\begin{rem}\label{r:irredn2}
The type $\tau_4$ is a full $\dual \F_q^\times$-orbit of size $q-1$ (under the action from Section \ref{s:cliff}), 
so the restriction to $\Sl_2(q)$ of every character of this type will remain irreducible according to Remark \ref{r:bijclif}.
\end{rem}

\begin{rem}\label{r:splittau1}
The action on a representation $R^G_T(\alpha,\beta)$ of type $\tau_1$ identifies the unordered 
pairs $\{\alpha, \beta\}$ with all pairs $\{\gamma\alpha, \gamma\beta\}$ for $\gamma\in\dual \F_q^\times,$ 
in particular with $\{\alpha/\beta, 1\}$. To simplify, we may assume the character is $R^G_T(\alpha,1),$ 
for a nontrivial $\alpha$. The orbit of $\{\alpha,1\}$ is made up of all the unordered pairs $\{\gamma\alpha,\gamma\}.$ 

Its stabilizer is trivial unless $\{\gamma\alpha,\gamma\}=\{\alpha,1\}$ for some $\gamma\neq 1$. That can only happen for $\gamma=\alpha^{-1}$ and 
$\alpha^2=1$. Let us call 
$\alpha_0\in\dual \F_q^\times$
the unique order two character. There is only one since the group $\dual \F_q^\times\=  \F_q^\times$ is cyclic. 
\end{rem}

\begin{rem}
We must point out that since both $\{\alpha,1\}$ and $\{1,\alpha^{-1}\}$ lie in the same orbit, $R^G_T(\alpha,1)$ 
and $R^G_T(\alpha^{-1},1)$ restrict to the same representation of $\Sl_2(q)$, so we only have $\dfrac{q-3}{2}$ of this type ($\alpha^2\ne 1$).
\end{rem}

In conclusion:
\begin{prop}
Every character of type $\tau_2$ restricts to an irreducible character of $\Sl_2(q)$, 
except that of $R^G_T(\alpha_0,1)$ whose restriction has two irreducible components. We will call them $\chi^+_{\alpha_0}$ and $\chi^-_{\alpha_0}.$
\end{prop}


\begin{rem}\label{r:splittau2}
With $-R_{T^s}^G(\omega)$ the situation is similar. The orbit of a representation corresponding to 
$\omega$ is given by those representations associated to the $\gamma\omega\in\Gamma_2$ 
where $\gamma \in \Gamma_1$ is regarded as in $\Gamma_2$  via the inclusion from \eqref{e:NGam}
$${}^TN_1^2:\Gamma_1 \hookrightarrow \Gamma_2.$$

Therefore, $(\omega\gamma)(a)$ should be computed as $\omega(a)\gamma(N_1^2(a))$ (which agrees with $\omega(a)\gamma(a^{q+1})=\omega(a)\gamma(a^{2})=\omega(a)\gamma^{2}(a)$ 
for $a\in\F_q$ in accordance with Remark \ref{r:evalua}). 

Since $\dual \F_q^\times$ acts by multiplication, we have $q-1$ different characters $\gamma\omega$. 
But we the irreducible representation $-R_{T^s}^G(\omega)$ depends only on the Frobenius orbit 
$\{\omega, \omega^q\},$ hence for a fixed $\omega\in\Gamma_2$ the stabilizer of the $\dual \F_q^\times$-action 
on the representations $-R_{T^s}^G(\gamma\omega)$ is given by the trivial character and all the $\gamma\in\Gamma_1$ 
such that $\omega\gamma=\omega^q$. And this can only happen when $\omega\gamma^2=\omega$ (since $\omega^{q^2}=\omega$ and $\gamma^q=\gamma$). 
Therefore, for $-R_{T^s}^G(\omega)$ to be stable by a nontrivial $\gamma$ we must have $\gamma=\alpha_0$ 
(the character of order $2$ from RemarK \ref{r:splittau1}) and $\omega^2\in \Gamma_1$ by Remark \ref{r:identify} and 
$$(\omega^2)^q=(\omega^q)^2=(\omega\alpha_0)^2=\omega^2.$$ 
 \end{rem}

Therefore, we conclude the following
\begin{prop}\label{p:splitting}
The representation $-R_{T^s}^G(\omega)$ will restrict to an irreducible representation of $\Sl_2(q)$ 
unless $\omega$ is the square root of a (non-square) degree one character, in which case its 
restriction have two components. 
\end{prop}

\begin{rem}
The map $\omega \mapsto \omega^{q-1}$ identifies characters $\omega,\gamma\omega$ in the same orbit and its kernel is $\Gamma_1$. Its image is the subgroup of $\Gamma_2$ of order $q+1.$ By Remark \ref{r:subquotdual} we can regard it as the dual of a cyclic group of order $q+1$ 
which we note $\mu_{q+1}.$ 
\end{rem}

\begin{rem}
The characters $\omega_0$ whose square have degree one will be mapped to those characters of order $1$ or $2$.
\end{rem}

By an abuse of notation we will keep calling them $\omega,\omega_0\in \mu_{q+1}$, $-R_{T^s}^G(\omega)$ is the restriction to $\Sl_2(q)$ of any representation $-R_{T^s}^G(\gamma\omega)$ of $\Gl_2(q)$. 
These are still irreducible for $\omega^2\ne 0$, and for the nontrivial $\omega_0\in\mu_{q+1}$ with $\omega_0^2=1$ we have 
two components and note its characters $\chi^+_{\omega_0}$ and $\chi^-_{\omega_0}.$

\begin{rem}
As before, since $-R_{T^s}^G(\omega)=-R_{T^s}^G(\omega^q)$  in the $\Gl_2$ case, we are going to have $-R_{T^s}^G(\omega)=-R_{T^s}^G(\omega^{-1})$ in the $\Sl_2$ case, so we only have $\dfrac{q-1}{2}$ of this type ($\omega^2\ne 1$).
\end{rem}

\newpage
Hence, the character table of $\Sl_2(q)$ looks like (see \cite{Fran}):
\noi
\begin{table}[h] \label{t:two}
\begin{center}
\caption{Character Table of $\Sl_2(q)$.}
\small
\begin{tabular}{|c|cccc|}
\hline
 classes & \mat a 0 0 a & \mat a b 0 a & \mat a 0 0 {a^{-1}} & \mat x 0 0 {x^q} \\
  &$ a=\pm 1$ & $a= \pm 1$ & $a\not\in \{1,-1\}$ & $x \neq  {x^q}$ \\
\hline
\# of classes & $2$ & $4$ & $ (q-3)/2 $&$ (q-1)/2 $\\
 class size & $1$ & $(q^2 - 1)/2$ & $q(q+1) $&$ q(q-1) $ \\
\hline
 $R_{T}^G(\alpha)$ & $(q+1)\alpha(a)$ & $\alpha(a)$ & $\alpha(a)+ \alpha(a^{-1})$ &$ 0 $\\
$\chi_{\alpha_0}^\pm$ & $\dfrac{(q+1)}{2}\alpha_0(a)$ & $\dfrac{\alpha_0(a)}{2} (1\pm \delanta)$ & $\alpha_0(a)$ &$ 0 $\\
$-R_{T^s}^G(\omega)$ & $(q-1)\omega(a)$ & $-\omega(a)$ & $0$ & $-(\omega(x)+\omega(x^q))$ \\
$\chi_{\omega_0}^\pm$ & $\dfrac{(q-1)}{2}\omega_0(a)$ & $\dfrac{\omega_0(a)}{2}(-1\pm \delanta)$ & $0$  & $-\omega_0(x)$ \\ 
 
 $1$ & $1$ & $1$ & $1$ & $1$ \\
$St_G$ & $q$ & $0$ & $1$ & $-1$ \\
\hline
\end{tabular}
\end{center}
\end{table}
\index{commands!environments!table}%

\noi
where $b$ is either $1$ or $y$ for some fixed non square $y\in\F_q\setminus (\F_q)^2$, $a\in\F_q^\times$ and 
$x\in\F_{q^2}^\times$ is a norm $1$ element of degree $2$ (in other words: $x^{q+1}=1$ and $x\neq \pm 1$), and 
in the $-R_{T^s}^G(\omega)$ row, the $\omega$ is taken in $\mu_{q+1}$ with $\omega^2\ne 1$.

The term $\delanta$ is a short for $\alpha_0(ab)\sqrt{\alpha_0(-1)q}$ and for our purposes is not 
really important since we only want the central values of the characters (i.e.: namely those of the first row).

\section{The number of $\F_q$-points of $\MM(\Sl_2)$}
\index{The number of $\F_q$-points of $\MM(\Sl_2)$
@\emph{The number of $\F_q$-points of $\MM(\Sl_2)$}}%

Thanks to formula \eqref{f:counting} we can compute the number $\NSol_2(q)=\NSMF_2(q)|H|$ \label{n:nsol}
of solutions 
of $$[y_1,x_1]\ldots [y_g,x_g]=z$$ for $x_i,y_i \in \Sl_2(q)$ and $z=-Id$ as:
\begin{eqnarray}\label{f:charsum2}
\NSol_2(q) &=& \sum_{\chi\in\Irr(H)}\chi(z)\left(\frac{|H|}{\chi(1)} \right)^{2g-1}
\end{eqnarray}
where $H=\Sl_2(q).$

By plugging the values from the first column of table \ref{t:two} in \eqref{f:charsum2}, 
and regarding $z$ as \mat{a}{0}{0}{a} for $a=-1$, the character sum 
\eqref{f:charsum2} becomes:
\begin{eqnarray}
\NSol_2(q) & = & \sum_{i=1}^{\frac{q-3}{2}} (q+1)(-1)^i\left(\frac{|H|}{q+1} \right)^{2g-1} \nonumber \\
& & +2\frac{q+1}{2}(-1)^\frac{q-1}{2}\left(\frac{|H|}{(q+1)/2} \right)^{2g-1} \nonumber \\
& & +\sum_{i=1}^{\frac{q-1}{2}} (q-1)(-1)^i\left(\frac{|H|}{q-1} \right)^{2g-1} \label{e:firstn2}\\
& & +2\frac{q-1}{2}(-1)^\frac{q+1}{2}\left(\frac{|H|}{(q-1)/2} \right)^{2g-1} \nonumber \\
& & +1\left(\frac{|H|}{1} \right)^{2g-1} \nonumber \\
& & +q\left(\frac{|H|}{q} \right)^{2g-1} \nonumber 
\end{eqnarray}
where different lines match the corresponding rows in the character table, for instance the second line is
$$\chi_{\alpha_0}^+(z)\left(\frac{|H|}{\chi_{\alpha_0}^+(1)} \right)^{2g-1} + \chi_{\alpha_0}^-(z)\left(\frac{|H|}{\chi_{\alpha_0}^-(1)} \right)^{2g-1}.$$

Dividing equation \eqref{e:firstn2} by $|H|$ and simplifying the sums we get 
\begin{eqnarray}
\NSMF_2(q) & = & 
 -(1+(-1)+\ldots +(-1)^{\frac{q-3}{2}-1})\pare{\frac{|H|}{q+1}}^{2(g-1)} \nonumber \\
& & +(-1)^\frac{q-1}{2} 2^{2g-1}\pare{\frac{|H|}{q+1}}^{2(g-1)} \nonumber\\
& & -(1+(-1)+\ldots +(-1)^{\frac{q-1}{2}-1}) \pare{\frac{|H|}{q-1}}^{2(g-1)} \label{e:secondn2}\\
& & +(-1)^\frac{q+1}{2}2^{2g-1} \pare{\frac{|H|}{q-1}}^{2(g-1)} \nonumber\\
& & +\pare{\frac{|H|}{1}}^{2(g-1)} +\pare{\frac{|H|}{q}}^{2(g-1)} \nonumber
\end{eqnarray}
which in turn reduces to:
\begin{eqnarray}
\NSMF_2(q) & = & 
 -\pare{\frac{1+(-1)^{\frac{q-1}{2}}}{2} } \pare{\frac{|H|}{q+1}}^{2(g-1)}\nonumber\\
& & +(-1)^\frac{q-1}{2}2^{2g-1}\pare{\frac{|H|}{q+1}}^{2(g-1)}\nonumber\\
& & -\pare{\frac{1+(-1)^{\frac{q+1}{2}}}{2}} \pare{\frac{|H|}{q-1}}^{2(g-1)}\label{e:quasipoly}\\
& & +(-1)^\frac{q+1}{2}2^{2g-1} \pare{\frac{|H|}{q-1}}^{2(g-1)}\nonumber\\
& & +\pare{\frac{|H|}{1}}^{2(g-1)} +\pare{\frac{|H|}{q}}^{2(g-1)}\nonumber
\end{eqnarray}
and finally to the quasi-polynomial:
\begin{eqnarray}
\NSMF_2(q) & = & 
 \left((-1)^\frac{q-1}{2}2^{2g-1}-\left(\frac{1+(-1)^{\frac{q-1}{2}}}{2}\right)\right)\left(\frac{|H|}{q+1} \right)^{2g-2} \nonumber\\
& & + \left((-1)^\frac{q+1}{2}2^{2g-1}-\left(\frac{1+(-1)^{\frac{q+1}{2}}}{2}\right)\right)\left(\frac{|H|}{q-1} \right)^{2g-2} \label{quasi}\\
& & +\pare{\frac{|H|}{1}}^{2g-2} +\pare{\frac{|H|}{q}}^{2g-2} \nonumber
\end{eqnarray}
which, for $q\equiv 1 \mod(4)$ becomes:
\begin{eqnarray}
\NSMF_2(q) & = &  
 \left(2^{2g-1}-1\right)\left(\frac{|H|}{q+1} \right)^{2g-2} - 2^{2g-1}\left(\frac{|H|}{q-1} \right)^{2g-2} \label{epoly}\\
& & +\pare{\frac{|H|}{1}}^{2g-2} +\pare{\frac{|H|}{q}}^{2g-2}. \nonumber
\end{eqnarray}
Since $|H|=q(q-1)(q+1)$ the number of points $\NSMF_2(q)$ depends polynomially in $q$ (provided $q\equiv 1 \mod(4)$),
with that polynomial having integer coefficients.

\section{Remarks}
\index{Remarks
@\emph{Remarks}}%

From the formulas obtained in last section we can deduce the following facts about ${\MM}(\Sl_2))$ and their $E$-polynomials:

\begin{enumerate}
\item Formula \eqref{epoly} implies that the Character Varieties ${\MM}(\Sl_2)$ have polynomial count and
thanks to Theorem \ref{t:Katz}, the polynomial expression found for $\NSMF_2$ are the $E$-polynomials 
${E(q;{\MM}(\Sl_2))}$ of $\MM(\Sl_2).$
\item 
For $g=1$, the quasi-polynomial computed in \eqref{quasi} tells us the number of pairs of matrices (up to $\PGl_2$-conjugation) 
in $\Sl_2(q)$ whose commutator is $-Id$ is $1$ because
\begin{eqnarray*}
\NSMF_2(q) & = &  
 (-1)^\frac{q-1}{2}2-\left(\frac{1+(-1)^{\frac{q-1}{2}}}{2}\right)\\
& & + (-1)^\frac{q+1}{2}2-\left(\frac{1+(-1)^{\frac{q+1}{2}}}{2}\right)\\
& & + 1 + 1 
\end{eqnarray*}
and since $\left((-1)^\frac{q-1}{2} + (-1)^\frac{q+1}{2}\right)=0,$
\begin{eqnarray*}
\NSMF_2(q) & = &  
 \left((-1)^\frac{q-1}{2} + (-1)^\frac{q+1}{2}\right)2 \\
 & & -\left(\frac{2+(-1)^{\frac{q-1}{2}}+(-1)^{\frac{q+1}{2}}}{2}\right)+ 2 \\ 
\NSMF_2(q) & = &  1.
\end{eqnarray*}

Therefore, the equation $[Y,X]=-Id$ has a
unique solution in $\Sl_2$ (modulo $\PGl_2$-action).


\item For $g\geq 2$ we see that the principal coefficient of the $E$-polynomial found in \eqref{epoly} 
comes from the term $\left(\frac{|H|}{1} \right)^{2g-2}$ and it is $1$.
This together with the conclusion of last remark (for the case $g=1$)
proves that the counting functions of the ${\MM}(\Sl_2)$ are monic of degree 
$\deg_q P = 3(2g-2).$ Then, the character varieties ${\MM}(\Sl_2)$ are connected of dimension $\dim(X).$ 
This and the next remark are particular cases of Corollary \ref{c:main}.

\item These $E$-polynomials are easily seen to be palindromic (i.e.: they remain unchanged when their coefficients are reversed). 
Being palindromic for a polynomial $P \in \Q[q]$ means $P(q)=q^{\deg P}P(q^{-1}).$ In our case, the polynomials $\NSMF_2$ 
have degrees $3(2g-2)$ and verify:
\begin{eqnarray*}
\NSMF_2(q^{-1})& = & 
 \left(2^{2g-1}-1\right)\left(q^{-1}(q^{-1}-1)\right)^{2g-2} 
 - 2^{2g-1}\left(q^{-1}(q^{-1}+1) \right)^{2g-2}\\
& & +\left(q^{-1}(q^{-1}+1)(q^{-1}-1)\right)^{2g-2} +\left((q^{-1}+1)(q^{-1}-1)\right)^{2g-2}
\end{eqnarray*}
multiplying both sides by $q^{6g-6}$
\begin{eqnarray*}
q^{6g-6}\NSMF_2(q^{-1})& = & 
 \left(2^{2g-1}-1\right)\left(q^{-1}(q^{-1}-1)q^{3}\right)^{2g-2} - 2^{2g-1}\left(q^{-1}(q^{-1}+1)q^{3} \right)^{2g-2}\\
& & +\left(q^{-1}(q^{-1}+1)(q^{-1}-1)q^{3}\right)^{2g-2} +\left((q^{-1}+1)(q^{-1}-1)q^{3}\right)^{2g-2} \\
q^{6g-6}\NSMF_2(q^{-1})& = & 
 \left(2^{2g-1}-1\right)\left(q(q-1)\right)^{2g-2} - 2^{2g-1}\left(q(q+1) \right)^{2g-2}\\
& & +\left((q+1)(q-1)\right)^{2g-2} +\left(q(q+1)(q-1)\right)^{2g-2}
\end{eqnarray*}
whose right-hand side is $\NSMF_2(q).$
\item The Euler Characteristic of $\MM(\Sl_2)$ is calculated evaluating at $q=1$ the polynomial $\NSMF_2(q)={E(q; {\MM}(\Sl_2))}$ from 
\eqref{epoly}, which in this case is
\begin{eqnarray*}\label{n:eval}
\NSMF_2(1) & = & 
 \left(2^{2g-1}-1\right)\pare{ q(q-1)}^{2g-2}\Bigg| _{q=1} - 2^{2g-1} \pare{ q(q+1)}^{2g-2}\Bigg| _{q=1}\\
& & +\pare{q(q-1)(q+1)}^{2g-2}\Bigg| _{q=1} +\pare{(q-1)(q+1)}^{2g-2} \Bigg| _{q=1}\\
 & = &   - 2^{4g-3}.
\end{eqnarray*}
We point out that the only term that survives (when evaluated at $q=1$) is the one corresponding to the characters of type $\tau_2$ that 
split when restricted to $\Sl_2(q).$

\end{enumerate}

\chapter{$E$-polynomial of $\MM(\Sl_n)$}\label{ch:epoly}
\index{$E$-polynomial of $\MM(\Sl_n)$%
@\emph{$E$-polynomial of $\MM(\Sl_n)$}}%

In this chapter we outline the computation of $E(q; {\MM}(\Sl_n)),$ the $E$-polynomials of the twisted 
character varieties $\MM(\Sl_n),$ by finding the counting functions $\NSMF_n(q)$ and proving they are polynomials
in $q$ for those prime powers $q$ satisfying \eqref{oddity}. 
Those polynomials turn out to be the sought $E$-polynomials thanks to Theorem \ref{t:Katz}.

The strategy to calculate $\NSMF_n(q)$ is the following
\begin{enumerate}
\item Write the character sum from \eqref{f:counting} as a sum in terms of $\Irr(\Gl_n(q))$.
\item Gather together the summands with the same type.
\item In last summation, separate the terms by splittings of the $\chi\in\Irr(\Gl_n(q))$ when restricted to $\Sl_n(q),$ according
to their stabilizers.
\item Filter the summands by degree of the classes of the characters $\gamma\in\Gamma$ from \ref{t:bijMac}. 
\item Reduce the sum to one in terms of characters $\alpha\in\Gamma.$
\item Go all the way back to the original sum with the results obtained.
\end{enumerate}
Throughout this chapter $g$ and  $n$ are positive integers, $q$ is a prime power with residue $1$ modulo $n,$
$\delta$ is a fixed generator of the cyclic group $ \Gamma_1=\wh{\F}_q^\times= \bra \delta,$ \label{n:deltagen}
$\Gl_n(q)$ will be noted by $G,$ 
$\Sl_n(q)$ will be noted by $H,$ 
$z$ is $\zeta_n\Id$ (an order $n$ central element in $H$), 
$\chi$ will note an irreducible character of $G,$  
$\theta$ an irreducible character of $H,$  
$s$ and $t$ will be sizes of orbits of associated characters verifying $st=q-1.$

We leave many of the calculations for Chapter \ref{ch:Kitchen}.

\section{First Reduction: From $\Sl_n$ to $\Gl_n$}
\index{First Reduction: From $\Sl_n$ to $\Gl_n$%
@\emph{First Reduction: From $\Sl_n$ to $\Gl_n$}}%

In order to apply Katz's result (Theorem \ref{t:Katz} here, or \cite{VilHaus}, appendix), 
we ought to compute $\NSMF_n(q),$ the number of $\F_q$-points of $\MM(\Sl_n).$ We will see that under condition
\eqref{oddity} this number depends polynomially in $q.$

By Corollary \eqref{c:sizesorbits}, Remark \ref{r:numberpoints}
and formula \eqref{f:counting} the number $\NSMF_n(q)$ of 
points of $\MM(\Sl_n)(q)$ is computed by the formula:
\begin{equation} \label{e:NSM}
\NSMF_n(q)=\frac{1}{|H|}\sum_{\theta\in \Irr(H)} \pare{ \frac{|H|}{\theta(1)} }^{2g-1}\theta(z).
\end{equation}

Thanks to the characterization of $\Irr(H)$ given in Section \ref{s:cliff},
we can express this summation in terms of characters from $\Irr(G).$ 

According to Remark \ref{r:bijclif}, there is a bijective correspondence between orbits under ${\quottt{G}{H}}$-conjugation 
of $\Irr(H)$ and orbits under $\dual{\pare{\quottt{G}{H}}}$-multiplication of $\Irr(G),$ and such correspondence was 
given by restricting to $H$ a character $\chi\in\Irr(G)$ and taking its irreducible components 
$\theta\in \Irr(H).$ 

Grouping together ${\quottt{G}{H}}$-orbits of $\Irr(H)$ in the sum \eqref{e:NSM} we get
\begin{eqnarray}
\NSMF_n(q) & = & \frac{1}{|H|}\sum_{\theta\in \Irr(H)}\pare{ \frac{|H|}{\theta(1)} }^{2g-1}\theta(z)\nonumber\\
 & = & \frac{1}{|H|}\sum_{\theta\in \Irr(H)}\pare{ \frac{|H|}{\theta(1)} }^{2g-1}\theta(z)\nonumber\\
 & = & \frac{1}{|H|}\sum_{[\theta] \in \Irr(H)/(G/H)} \pare{ \frac{|H|}{[\theta](1)} }^{2g-1}\sum_{\theta \in[\theta]}\theta(z)\label{e:NSM2}
\end{eqnarray}
where $[\theta]\subseteq \Irr(H)$ stands for a $\quottt{G}{H}$-orbit and $[\theta](1)$ is the degree of any element in such orbit.

Passing to the $\Irr(G)$ side by means of the correspondence from Remark \ref{r:bijclif}, and using that:
\begin{itemize}
\item $\displaystyle{\sum_{\theta \in[\theta]}\theta(z)=\chi(z)} $ for $\chi\in\Irr(G)$ the corresponding character whose $H$ 
restriction affords $[\theta],$
\item The degree $\chi(1)$ is $t\;\theta(1),$ where 
\begin{equation}\label{e:t}
t=t(\chi):=\size{[\theta]}, 
\end{equation}
\end{itemize}
we end up with:
\begin{eqnarray}
\NSMF_n(q) & = & \frac{1}{|H|}\sum_{[\chi]/\chi\in \Irr(G)}\pare{ \frac{|H|}{\chi(1)} }^{2g-1}t(\chi)^{2g-1}\chi(z) \label{e:NSM3}
\end{eqnarray}
where $[\chi]$ is the orbit of $\chi\in\Irr(G)$ of some size $s=s(\chi):=|[\chi]|$, $t$ is the size of 
the ${\quottt{G}{H}}$-orbit $[\theta]$ corresponding to the constituents of the restriction $\chi_H$ and their 
sizes satisfy $st=q-1$ (since ${\quottt{G}{H}} = \F_q^\times$ is cyclic of order $q-1$).

\begin{rem}\label{r:tdivn}
A character $\chi\in\Irr(G)$ whose $\wh{\pare{\quottt{G}{H}}}$-orbit has size $s$ must necessarily have a stabilizer of order $t=(q-1)/s.$
Since $\wh{\pare{\quottt{G}{H}}}=\wh \F_q^\times = \Gamma_1$ is cyclic, such stabilizer is generated by $\delta^s$ (Remark \ref{r:group}).
By Remark \ref{r:action} we know $\delta^s$ acts on the support of the $\Ji$-partitions $\L:\Gamma\to\PP$ associated to $\chi.$
Since multiplication is a free action, by Theorem \ref{t:bijMac} we conclude that $t$ divides each product $dm_{d,\l},$ and 
hence divides $n$.
\end{rem}


\section{Second Reduction: From $\Gl_n$ to types}
\index{Second Reduction: From $\Gl_n$ to types%
@\emph{Second Reduction: From $\Gl_n$ to types}}%

We gather together from \eqref{e:NSM3} those characters $\chi\in\Irr(G)$ having the same type $\tau.$ 

Using formula \eqref{e:central},
equation \eqref{e:NSM3} becomes
\begin{eqnarray}
\NSMF_n(q) & = &\frac{1}{|H|}\sum_{\substack{[\chi],t \\ \chi\in \Irr(G)\\ t=t(\chi) }}\pare{ \frac{|H|}{\chi(1)} }^{2g-1}t^{2g-1}\chi(z)\nonumber\\ 
 & = & \sum_{\substack{[\chi],t \\ \chi\in \Irr(G)\\ t=t(\chi) }}\pare{ \frac{|H|}{\chi(1)} }^{2(g-1)}t^{2g-1}\Delta_\chi(\zeta_n)\nonumber\\ 
 & = & \sum_{\substack{\tau, t\\ |\tau|=n \\ t|n}}\pare{ \frac{|H|}{\chi_\tau(1)} }^{2(g-1)}t^{2g-1}
 \sum_{\substack{[\chi] \in \quot{\Irr(G)}{\dual{\F_q^\times}}\\ t(\chi)=t \\ \tau(\chi)=\tau }}\Delta_\chi(\zeta_n)\nonumber\\
 & = & \sum_{\substack{\tau, t\\ |\tau|=n \\ t|n}}\pare{ \frac{|H|}{\chi_\tau(1)} }^{2(g-1)}t^{2g-1}
 C_\tau^t \label{e:NSM4} 
\end{eqnarray}
where 
\begin{equation} \label{e:Ctaut}
C_\tau^t:=\sum_{\substack{[\chi] \in \quot{\Irr(G)}{\wh{\F_q^\times}}\\|[\chi]|=s\\ \tau(\chi)=\tau}}\Delta_\chi(\zeta_n) 
\end{equation}
and $\displaystyle \Delta_\chi(\zeta_n)$ is the $\Delta_\L(\zeta_n)=\prod_{\gamma}\gamma^{|\L(\gamma)|}(\zeta_n)$ (evaluated in $\Gamma_1$), from 
\eqref{e:central}.

Using that $|H|=(q-1)^{-1}|G|$ and plugging in \eqref{e:quot} in \eqref{e:NSM4} we find
\begin{equation} \label{e:NSM5}
\NSMF_n(q)=\sum_{\tau,t}\pare{q^{\frac{n^2}{2}}\frac{\H_{\tau'}(q)}{q-1} }^{2(g-1)}t^{2g-1}C_\tau^t.
\end{equation}

\begin{rem}\label{r:epolyCtaut}
In Section \ref{s:wrap} we will see that $C^t_\tau$ does not depend on $q$ (at least as soon as $q$ 
satisfies \eqref{oddity}) and hence it is just a number.
This, together with the fact that the factors
$\pare{ \frac{|H|}{\chi_\tau(1)} }=\pare{ q^{\frac{n^2}{2}}\frac{\H_{\tau'}(q)}{q-1}}$ are polynomials in $\Z[q]$ 
(formulas \eqref{e:hookpolyl}, \eqref{e:hookpolyL}) and Remark \ref{r:integercoeff}),
will lead us to the remarkable fact that $\NSMF_n(q)$ is a polynomial for those values of $q.$ Then, by Remark \ref{r:spreadout}, 
the Character Variety $\MM(\Sl_n)$ results polynomial count.
\end{rem}


\begin{rem}
For $g\geq 2$ the highest degree term in \eqref{e:NSM5} is the one corresponding 
to the type $\tau_p$ with unique nonzero multiplicity $(\tau_p)_{1,(1,\ldots,1)}=1,$ 
namely that of $1_G,$ the trivial character\footnote{also: the biggest $\pare{\frac{|H|}{\chi(1)}}$ is attained by the character with smallest degree $\chi(1).$} (see Remarks \ref{r:principal} and \ref{r:1diml}).
\end{rem}

\begin{rem}
Since the $\dual{\pare{\quottt{G}{H}}}$-orbit of $1_G$ is the set of all $\tau_p$-type 
characters in $\Irr(G)$, and has size $s=q-1$, the corresponding ${\quottt{G}{H}}$-orbit will have 
size $t=1$ (Remark \ref{r:bijclif}), hence their $H$-restriction are irreducible\footnote{the restriction of the trivial character $1_G$ is $1_H,$ hence it cannot have more than one component.}.
\end{rem}

\begin{rem}\label{r:degepoly}
The only $t$ for which $C^t_{\tau_p}\ne 0$ is $t=1.$ Since the sum \eqref{e:Ctaut} giving $C^1_{\tau_p}$ has only one summand, 
whose value is $1,$ the principal coefficient of $\NSMF_n(q)$ will be $1$ and its degree $(2g-2)(n^2-1).$ 
This, together with Corollary \ref{c:palin} proves Corollary \ref{c:main}.
\end{rem}

\begin{rem}\label{r:tautau}
From \eqref{e:Ctaut} and the fact that $\size{\L(\gamma)}=\size{\L'(\gamma)}$ we see that 
$C_\tau^t$ and $C_{\tau'}^t$ have the same value.
\end{rem}

\section{Third Reduction: Controlling Stabilizers}
\index{Third Reduction: Controlling Stabilizers%
@\emph{Third Reduction: Controlling Stabilizers}}%


\subsection{Fixed Stabilizer}
\index{Fixed Stabilizer%
@\emph{Fixed Stabilizer}}%

Let us focus now in the calculation of $C^t_\tau.$ Since they are defined in \eqref{e:Ctaut} by a sum 
on $\dual{\pare{\quottt{G}{H}}}$-orbits of characters, all of which have size $s=(q-1)/t$, we can transform such 
sum in one on characters by overcounting the character classes:
\begin{equation} \label{e:sCtaut}
sC_\tau^t=\sum_{\substack{\chi\in \Irr(G)\\|[\chi]|=s\\ \tau(\chi)=\tau}}\Delta_\chi(\zeta_n)
\end{equation}
where $\dual{\pare{\quottt{G}{H}}}$-conjugate characters are considered separately. The others restrictions on the sum remain the same, 
in other words, we are summing over characters $\chi\in\Irr(G)$ (instead of character classes) 
of type $\tau$ and whose stabilizers have order exactly $t$.


\begin{rem}
For a $\chi\in\Irr(G)$ to have stabilizer of order $t$ is the same as verifying $\chi=\delta^s\chi$ and $s$ the smaller such power of $\delta$ that stabilizes it.
\end{rem}

\subsection{Bounded Stabilizer}
\index{Bounded Stabilizer%
@\emph{Bounded Stabilizer}}%

It will be convenient to work with characters that are stable by a fixed power of $\delta,$ regardless 
of the actual size of its stabilizer. And that amounts to take a cumulative sum of these $sC^t_\tau$ which we call:

\begin{equation}\label{e:Shat}
\wh{S}(t,\tau):=\sum_{\substack{k|\frac{n}{t}\\ n=|\tau|}}\frac{s}{k}C_\tau^{tk}=\sum_{\substack{\chi\in \Irr(G)\\ \delta^s\chi=\chi\\ \tau(\chi)=\tau}}\Delta_\chi(\zeta_n). 
\end{equation}

\begin{rem}\label{r:ShattoCtaut}
Thanks to the M\"obius Inversion Formula, we can solve for $C^t_\tau$ in \eqref{e:Shat} getting
\begin{equation}\label{e:ShattoCtaut}
C_\tau^t=\frac{t}{q-1}\sum_{\substack{k|\frac{n}{t}\\ n=|\tau|}}\mu(k)\wh{S}(kt,\tau). 
\end{equation} 
\end{rem}

\begin{rem}\label{r:ShattoCo}
Plugging $t=1$ in \eqref{e:Shat}, and keeping in mind that $\delta^s$ becomes trivial in this case, we get
$$
\wh{S}(1,\tau)=\sum_{\substack{\chi\in \Irr(G)\\ \tau(\chi)=\tau}}\Delta_\chi(\zeta_n)=C_\tau
$$
were the coefficients $C_\tau$ are those form \cite{VilHaus} (3.2.2).
\end{rem}

\section{Fourth Reduction: The New-degree}
\index{Fourth Reduction: The New-degree%
@\emph{Fourth Reduction: The New-degree}}%

In order to evaluate the $\Delta_\chi$'s, we need to consider one extra feature we shall define next.

\begin{defn}\label{n:newdeg}
The Frobenius automorphism on $\Gamma$ induces an action in the quotient $Frob_q \acts \quott{\Gamma}{\bra{ \delta^s }}$.
Define the \emph{newdegree} of $\alpha\in\Gamma$ relative to $t$ as the degree of its class 
$[\alpha]\in \quott{\Gamma}{\bra {\delta^s} },$ namely
$$\deg^t(\alpha):=\deg([\alpha])=\size{\set{ [\alpha],[\alpha]^q,[\alpha]^{q^2},\ldots } }.$$
\end{defn}

%
%

\begin{rem}\label{r:tdhat}
Given $t|n$ and $\alpha$ of degree $d$ and newdegree $\dhat$ we have 
$$\dhat t=\size{ \pare{ \bra{Frob_q}\times\bra{\delta^{s}}}\alpha}=\size{\set{\delta^{si}\alpha^{q^j}}_{i,j}}=d\that$$ 
for some integer $\that.$
\end{rem}

\begin{rem}\label{r:tdhatdivn}
Given a $\chi\in\Irr(G)$ stable by multiplication by $\delta^s$ and $\L$ the corresponding
multipartition, we have that the product $\dhat t$ from last remark divides $dm_{d,\l}(\L)$ for any $\l \in \PP,$ since 
for a fixed $\alpha\in\Gamma$ every $i$ and $j,$ the partitions $\L(\delta^{si}\alpha^{q^j})$ are the same. 
\end{rem}

\begin{rem}\label{r:dhatdividesd}
As a particular case of Lemma \ref{l:grthy} from Appendix \ref{ap:twolem} we get that $\dhat$ divides $d$, and hence, by last 
remark $\that$ divides $t$. This follows by realizing that $d$ is the size of a $Frob_q$-orbit 
of a character, and that $\dhat$ is the size of a $Frob_q$-orbit of character classes.
\end{rem}

\begin{defn}
We say that a type $\tau$ is pure of degree $d$ if all of its 
multiplicities $(\tau)_{d',\l}$ are $0$ for $d'\neq d.$ 
This is the same as saying that all the $\gamma\in\Gamma$ in the support of a $\L$ of type $\tau$ have degree $d.$ 

We also say a character $\chi$ is pure if its type is so. 
\end{defn}

\begin{defn}
Let us consider $t$ some divisor of $n$ and $\chi$ a character of $G$ corresponding to a multi-partition $\L\in\PP_n^{\Ji}.$ 
Let us further assume $\chi$ stable by $\delta^s$\footnote{remember $st=q-1$}. We say 
$\chi$ is pure of newdegree $\dhat$ if all the $\gamma$'s in the support of $\L$ have newdegree $\dhat.$
\end{defn}
\begin{notation}
In this case we write $\deg^t(\chi)=\dhat.$ 
\end{notation}

\begin{rem}\label{r:claimShat}
We claim $\wh{S}(t,\tau)=0$ unless $\tau$ is pure of some degree $d$, in particular $d|n$. This will be proven in Lemma \ref{l:whS} 
of next Chapter.
\end{rem}

Writing
\begin{equation}\label{e:Sdef}
S(t,\dhat,\tau):=\sum_{\substack{\chi\in \Irr(G)\\ \delta^s\chi=\chi\\ \deg^t(\chi)=\dhat\\ \tau(\chi)=\tau}}\Delta_\chi(\zeta_n)
\end{equation}
we get
\begin{equation}\label{e:Sform}
\wh{S}(t,\tau)=\sum_{\dhat|d}S(t,\dhat,\tau) 
\end{equation}
where $d$ is the degree of all the characters in the support of $\tau$. Here we are implicitly using the fact that the 
characters of mixed newdegrees add up to zero, which will be proved in Lemma \ref{l:mixnewdeg}.


\section{Fifth Reduction: Writing it in terms of $\Gamma$}
\index{Fifth Reduction: Writing it in terms of $\Gamma$%
@\emph{Fifth Reduction: Writing it in terms of $\Gamma$}}%

\subsection{Ordering the characters}
\index{Ordering the characters%
@\emph{Ordering the characters}}%

Now take some character $\chi$ corresponding to some function $\L\in\PP^\Ji_n$, pure of degree $d$ and of 
newdegree $\dhat$, for some given $t$ divisor of $n$.

Recalling \eqref{e:deltafactorcentral}
\begin{eqnarray}
\Delta_\L(\zeta_n) & = & \prod_{\{\gamma\}\in \Ji}\gamma^{|\L(\gamma)|}(\zeta_n) \nonumber
\end{eqnarray}
the product being taken over $Frob_q$-orbits, and the evaluations are done in $\Gamma_d.$

Writing it in terms of $\dual \F_q^\times$-orbits we have
\begin{eqnarray}
\Delta_\L(\zeta_n) & = & \prod_{\{[\gamma]\}\in \quott{\Gamma}{\langle Frob_q, \delta^s\rangle } } 
\prod_{i=0}^{\that-1} (\delta^{si}\gamma)^{|\L(\gamma)|}(\zeta_n) \nonumber\\
& = & 
\prod_{\{[\gamma]\}\in \quott{\Gamma}{\langle Frob_q, \delta^s\rangle}} 
\pare{\gamma^{\that|\L(\gamma)|}(\zeta_n) \prod_{i=0}^{\that-1} (\delta^{dsi})^{|\L(\gamma)|}(\zeta_n) } \nonumber
\end{eqnarray}
where the $\delta$ is evaluated in $\Gamma_1,$ and that is why we add the factor $d$ to the exponent.
Grouping together and pulling out all the powers of $\delta$ we get
\begin{eqnarray}
\Delta_\L(\zeta_n) & = & 
\prod_{\{[\gamma]\}\in \quott{\Gamma}{\langle Frob_q, \delta^s\rangle}}\gamma^{|\L(\gamma)|}(\zeta_n^{\that}) \delta(\zeta_n)^{ds{|\L(\gamma)|}\sum_{i=0}^{\that-1}i} \nonumber\\
& = & 
\prod_{\{[\gamma]\}\in \quott{\Gamma}{\langle Frob_q, \delta^s\rangle}}
\gamma^{|\L(\gamma)|}(\zeta_n^{\that}) \delta(\zeta_n)^{ds{|\L(\gamma)|}{\that\choose 2}} \nonumber\\
& = & 
\delta(\zeta_n)^{ds{\that\choose 2}{\sum|\L(\gamma)|}}\prod_{\{[\gamma]\}\in \quott{\Gamma}{\langle Frob_q, \delta^s\rangle}}\gamma^{|\L(\gamma)|}(\zeta_n^{\that}) \label{e:deltafactor}
\end{eqnarray}
where the summation $ds{\that\choose 2}{\sum|\L(\gamma)|}$ in the first exponent has the same 
indexing as the product, namely the $Frob_q$-orbits of the classes of characters 
$[\gamma]\in\quott{\Gamma}{\bra{\delta^s}}.$ The $\delta$ is evaluated as in $\Gamma_1$ and the $\gamma$'s are evaluated
as in $\Gamma_d.$

We next overcount every $\chi$ in the summation \eqref{e:Sdef} by ordering the $\gamma$'s in tuples of character classes 
$[\alpha_i]\in\quott{\Gamma}{\bra{\delta^s}}$ and distinguishing between $Frob_q$-conjugates as well and then we express:
\begin{equation}\label{e:ZtoS}
S(t,\dhat,\tau)=\frac{\delta(\zeta_n)^{ds{\that\choose 2}{\sum \l_i}}}{\dhat^m\prod \wh{m_\l}!}
Z(d,t,\dhat,\vec{\l})
\end{equation}
where $\wh{m_\l}$ is defined by $\dhat t \wh{m_\l}= m_{d,\l}d$, $\vec{\l}=(\l_1,\ldots,\l_m)$ is 
the list of the exponents corresponding to a fixed ordering of the $|\L(\gamma)|$'s, and
\begin{equation}\label{e:Zdef}
Z(d,t,\dhat,\vec{\l}):=\sum_{\substack{([\alpha_1],\ldots,[\alpha_m])\in (\quott{\Gamma_d}{\langle\delta^s\rangle})^m\\ \{[\alpha_i]\}\neq\{[\alpha_j]\},\: i\neq j \\ \deg(\alpha_i)=d \\ \deg([\alpha_i])=\dhat}}\prod_{i=1}^m\alpha_i^{\l_i}(\zeta_n^{\that})
\end{equation}
where the $\alpha_i$'s are evaluated as in $\Gamma_d.$

\begin{rem}
The range of the sum should be read as follows:
\noi
\begin{itemize}
\item We are taking $m$-tuples of classes of characters $[\alpha]$ were $\alpha\in\Gamma_d.$ We identify two of 
them whenever their entries differ only by factors of $\delta^s$, namely if they are $\delta^s$-conjugate.
\item Different entries in the $m$-tuple cannot be $Frob_q$-conjugates. In particular they correspond to 
different $\alpha_i\in\Gamma_d.$
\item All the entries have degree exactly $d.$
\item All the entries have newdegree exactly $\dhat$ relative to $t.$ That is to say, the degree of the 
class $[\alpha_i] \in \quott{\Gamma}{\bra{\delta^s}}$ is $\dhat.$
\end{itemize}
\end{rem}

\begin{rem}\label{r:deltafactor}
The factor $\delta(\zeta_n)^{ds{\that\choose 2}{\sum|\L (\gamma)|}}$ is $\pm 1$ and it is $1$ under 
condition \eqref{oddity}
.
This will say that for those prime powers $q\equiv 1 \mod(n)$ (i.e.: those for which $\MM(\Sl_n)$ is defined) 
the function $\NSMF_n(q)$ is a quasi-polynomial in $q$, with period at most $2.$
\end{rem}

\begin{rem}\label{r:evawell}
Since the evaluations $\alpha_i(\zeta_n^{\that})$ are done in $\Gamma_d$, they are well defined in $\quott{\Gamma}{\bra{\delta^s}}$ because
$\zeta_n^{\that}\in\Ker({}^TN_1^d(\delta^s))$ (see Remark \ref{r:subquotdual}).
This last assertion holds since
$${}^TN_1^d(\delta^s)(\zeta_n^{\that})=\delta^s(N_1^d(\zeta_n^{\that}))=\delta^s(\zeta_n^{d\that})=\delta^{sd\that}(\zeta_n)$$
evaluated as in $\Gamma_1$, which is $1$ as $sd\that=st\dhat=(q-1)\dhat$ is divisible by $n.$ 
\end{rem}


\subsection{Ignoring degrees}
\index{Ignoring degrees%
@\emph{Ignoring degrees}}%

Now define the accumulated sums of these $Z(d,t,\dhat,\vec{\l})$ in order to get rid of the 
condition on the degrees of the characters $\alpha_i.$ Note that this will not remove the condition on 
their newdegrees.

\begin{defn}
Keeping in mind that the degree of a character must be divisible by its newdegree $\dhat$ (Remark \ref{r:dhatdividesd}) we define:
\begin{equation}\label{e:Zhatd}
\wh{Z}(D,t,\dhat,\vec{\l}):=\sum_{\substack{d|D\\\dhat|d }}Z(d,t,\dhat,\vec{\l}).
\end{equation}
\end{defn}

\begin{rem}\label{r:ZhattoZ}
Then, inverting \eqref{e:Zhatd} with formula \eqref{f:mIdiv} we will have:
\begin{equation}\label{e:ZhattoZ}
Z(D,t,\dhat,\vec{\l})=\sum_{\substack{d|D\\ \dhat|d }}\mu\pare{\frac{D}{d}}\wh{Z}(d,t,\dhat,\vec{\l}). 
\end{equation} 
\end{rem}

\begin{rem}
We need to point out that each term $Z(d,t,\dhat,\vec{\l})$ corresponds to a sum over $m$-tuples of classes of characters of $\Gamma_d$ 
where the $d$ changes, and hence so does the exponent $\that$, being $\dhat t / d.$ These two changes compensate 
each other in the sense that for a given $\alpha \in \Gamma_d$ evaluating 
at $\zeta_n^{\dhat t / d}$ regarding $\alpha$ as in $\Gamma_d$ and evaluating 
at $\zeta_n^{\dhat t / D}$ regarding $\alpha$ as in $\Gamma_D$ yields to the same result since 
in the second case one needs to compose with $N_d^D$ which amounts to multiplying the exponent of $\zeta_n$ by a factor of $D/d.$
\end{rem}

From this remark we conclude that the sum defining $\wh{Z}(D,t,\dhat,\vec{\l})$ could be computed with 
a fixed $\that=\dfrac{t\dhat}{D}$ and all the evaluations with the $\alpha$'s regarded as in $\Gamma_D$ 
even when the $\alpha$ has a smaller degree.

Thus, we could compute $\wh{Z}$ as:
\begin{equation}\label{e:Zhatf}
\wh{Z}(D,t,\dhat,\vec{\l})=\sum_{\substack{([\alpha_1],\ldots,[\alpha_m])\in \pare{\quott{\Gamma_D}{\bra{\delta^s}}}^m\\ \{[\alpha_i]\}\neq\{[\alpha_j]\},\: i\neq j \\ \deg([\alpha_i])=\dhat}}\prod_{i=1}^m\alpha_i^{\l_i}(\zeta_n^{\that}). 
\end{equation}

\begin{rem}\label{r:Zhatf}
Note that in \ref{e:Zhatf} there are summands $\prod_{i=1}^m\alpha_i^{\l_i}(\zeta_n^{\that})$ 
not all of whose factors $\alpha_i$ have the same degree. It does not affect the 
value of the sum, because all this terms with mixed degrees will just cancel each other (Lemma \ref{l:mixdeg})
\end{rem}

\section{Wrapping up}\label{s:wrap}
\index{Wrapping up%
@\emph{Wrapping up}}%

In section \ref{s:values} we compute the values of the $\wh{Z}(D,t,\dhat,\vec{\l})$ and get the:

\begin{lem}\label{l:keylemma}
Given $D,$ $\dhat$ and $t$ divisors of $n$ and a list of exponents $\vec{\l}$ 
verifying  $ D | t \dhat,$ $\dhat|D$ and $n=\dhat t \sum\l_i,$ for $q$ satisfying \eqref{oddity} we have:
\begin{equation} \label{e:keylemma}
\wh{Z}(D,t,\dhat,\vec{\l})=\left\{
\begin{array}{cc}
\mu(\dhat)(-\dhat)^{m-1}(m-1)!\frac{q-1}{t} & \text{when } \gcd(D,t)=1  \\
0 & \text{otherwise}  \end{array}
\right.
\end{equation}
\end{lem}

\begin{rem}\label{r:dhatD}
Since $\wh{Z}(D,t,\dhat,\vec{\l})$ is zero unless the three conditions $ D | t \dhat $, $\dhat|D$
and $\gcd(D,t)=1$ hold, the only nonzero values of $\wh Z$ arise when $\dhat = D$. 
\end{rem}

Thus, under the condition \eqref{oddity} for $q$, and in view of last remark, formula \eqref{e:ZhattoZ} becomes:
\begin{eqnarray}
Z(D,t,\dhat,\vec{\l})&=&\mu\pare{\frac{D}{\dhat}}\wh{Z}(\dhat,t,\dhat,\vec{\l}) \nonumber \\
&=&\mu\pare{\frac{D}{\dhat}}\mu(\dhat)(-\dhat)^{m-1}(m-1)!\frac{q-1}{t} \label{f:Zexp}
\end{eqnarray}
provided $(\dhat;t)=1$, and zero otherwise.

Then, thanks to \eqref{f:Zexp} and Remark \ref{r:deltafactor}, formula \eqref{e:ZtoS} becomes:
\begin{equation}\label{f:Sexp}
S(t,\dhat,\tau)=\left\{
\begin{array}{cc}
\mu(\frac{d}{\dhat})\frac{q-1}{t}(-1)^{m-1}\frac{\mu(\dhat)}{\dhat} \frac{(m-1)!}{\prod \wh{m_\l}!} & \text{when } \gcd(\dhat,t)=1  \\
0 & \text{otherwise}  \end{array}
\right.
\end{equation}
where $d$ is the degree of $\tau$.

%
%
%
%
%
%
%
%

We recognize the factor 
\begin{equation}\label{f:factor}
(-1)^{m-1}\frac{\mu(\dhat)}{\dhat} \frac{(m-1)!}{\prod \wh{m_\l}!} 
\end{equation}
as $C_{\wh{\tau}}/(q-1)=\Co_{\wh{\tau}}$ (calculated in \cite{VilHaus} (3.4.2), see \eqref{f:Ctau} in Remark \ref{r:Ctau} below) 
for some other type ${\wh{\tau}},$ depending on $\tau$, $t$ and $\dhat$. 

Let's describe this new ``quotient'' type.
\begin{defn}
The type $\tau$ has only characters of degree $d$, and we are dividing them by a factor $t_d$ to get $\dhat$.
We are also dividing the multiplicities $m_{d,\l}$ in order to get the $\wh{m}_\l$ in such a way to have
\begin{eqnarray*}
d m_{d,\l} & = & \dhat t \wh{m}_\l  \\
\frac{t_d}{t} m_{d,\l} & = & \wh{m}_\l.
\end{eqnarray*}
\end{defn}

\begin{rem}
Writing $t_m$ for the factor $\frac{t}{t_d},$ we verify 
\begin{itemize}
\item $t=t_dt_m,$
\item the size of $\wh \tau$ is  $\displaystyle{\size{\wh \tau}=\dfrac{n}{t},}$
\item its degree $\deg(\wh \tau)$ is $\displaystyle{\dhat=\frac{d}{t_d}}$ and
\item the multiplicities $(\wh \tau)_{\dhat,\l}$ are $\displaystyle{\wh{m}_\l = \frac{(\tau)_{d,\l}}{t_m}.}$
\end{itemize}
\end{rem}

\begin{notation}\label{n:quotype}
We write $\tau/(t_d,t_m)$ for this new type $\wh \tau$.
\end{notation}
 
Formula \eqref{f:Sexp} then becomes
\begin{equation}\label{f:C1toS}
S(t,\dhat,\tau)=\left\{
\begin{array}{cc}
\frac{q-1}{t}\mu(t_d)\Co_{\tau/(t_d,t_m)} & \text{when } \gcd(\dhat,t)=1  \\
0 & \text{otherwise}  \end{array}
\right.
\end{equation}

And for $\tau$ pure of degree $d$, and $t$ a divisor of $n$, plugging \eqref{e:Sform} in \eqref{e:ShattoCtaut}, 
the coefficient $C_\tau^t$ becomes
\begin{eqnarray}
C_\tau^t & = & \frac{t}{q-1}\sum_{k|\frac{n}{t}}\mu(k)\wh{S}(kt,\tau)\nonumber\\
& = & \frac{t}{q-1}\sum_{k|\frac{n}{t}}\mu(k)\sum_{\dhat|d}S(kt,\dhat,\tau)\nonumber
\end{eqnarray}
which, thanks to \eqref{f:C1toS}, boils down to 
\begin{eqnarray}
C_\tau^t & = & \frac{t}{q-1}\sum_{k|\frac{n}{t}}\mu(k)\sum_{\substack{t_d,t_m\\ t_dt_m=tk\\ \exists {\tau/(t_d,t_m)}\\(\dhat,tk)=1} }
\frac{q-1}{tk}
\mu(t_d)\Co_{\tau/(t_d,t_m)}\nonumber\\
& = & \sum_{\substack{ t_d,t_m,k \\ \exists {\tau/(t_d,t_m)} \\ (\dhat,t_dt_m)=1\\ t_dt_m = tk}}
\frac{\mu(k)}{k}\mu(t_d)\Co_{\tau/(t_d,t_m)} \label{f:C1toCtaut}
\end{eqnarray}
where this last sum is taken over all quotient types $\tau/(t_d,t_m)$ of degree $\dhat$, 
such that $t_dt_m$ is multiple of $t$ and is relatively prime to $\dhat$.

\begin{rem}\label{r:Ctau}
From \ref{r:ShattoCo}, \eqref{e:Sform}  and \eqref{f:Sexp} we recover formula (3.4.2) from \cite{VilHaus},
\begin{equation}\label{f:Ctau}
C_\tau=\left\{
\begin{array}{cc}
(-1)^{m-1}({q-1})\frac{\mu(d)}{d} \frac{(m-1)!}{\prod {m_{d,\l}}!} & \text{for $\tau$ pure of degree $d$} \\
0 & \text{otherwise}  \end{array}
\right.
\end{equation}
\end{rem}


\section{Proof of Main Statements}\label{s:prfmthm}
\index{Proof of Main Statements%
@\emph{Proof of Main Statements}}%


We are now in conditions to prove Theorem \ref{t:main}.
\proof
The separated scheme from \eqref{d:spreadoutcharvar} gives us a spreading out of $\MM(\Sl_n)$ thanks to
\eqref{e:sesh}. 

By Remark \ref{r:epolyCtaut} and \eqref{f:Ctau}, the counting functions $\NSMF_n(q)$ are polynomials in $q.$

This, together with \eqref{e:kac} proves that the Character Varieties $\MM(\Sl_n)$ have polynomial count.
Finally, by Katz's Theorem (\ref{t:Katz}) the theorem follows.
\endproof

As an immediate consequence we have Corollary \ref{c:main} the proof of which we shall write now.

\proof
Plugging expression \eqref{f:C1toCtaut} in formula \eqref{e:NSM5} and writing $\wh{\tau}$ for 
the quotient types $\tau/(t_d,t_m)$ we get:
\begin{eqnarray}
\NSMF_n(q)& = & \sum_{\substack{\tau,t \\ t | n \\ |\tau|=n  }}\pare{\frac{q^{\frac{n^2}{2}}\H_{\tau'}(q) }{q-1} }^{2(g-1)}t^{2g-1}C_\tau^t \nonumber\\
%
& = & \sum_{\substack{\tau,t \\ t | n \\ |\tau|=n  }}\pare{\frac{q^{\frac{n^2}{2}}\H_{\tau'}(q) }{q-1} }^{2(g-1)}t^{2g-1}
\sum_{\substack{ t_d,t_m,k \\ \exists {\tau/(t_d,t_m)} \\ (\dhat,t_dt_m)=1\\ t_dt_m = tk}}
\frac{\mu(k)}{k}\mu(t_d)\Co_{\tau/(t_d,t_m)}\nonumber\\
& = & \sum_{\substack{\tau,t \\ t | n \\ |\tau|=n  }}
\sum_{\substack{ t_d,t_m,k,\wh \tau \\ \wh \tau = {\tau/(t_d,t_m)} \\ (\dhat,t_dt_m)=1\\ t_dt_m = tk}}
t^{2g}\frac{\mu(k)}{kt}\mu(t_d)\Co_{\wh{\tau} }\pare{\frac{q^{\frac{n^2}{2}}\H_{\tau'}(q) }{q-1} }^{2(g-1)}\nonumber
\end{eqnarray}
parameterizing the sum by pairs of types $\tau$, $\wh\tau$ one of which is a quotient of the other
\begin{eqnarray}
\NSMF_n(q)& = & \sum_{\substack{ \tau,\wh \tau ,t_d,t_m\\ |\tau|=n \\ \wh{\tau} = \tau/(t_d,t_m) \\ (\dhat,t_dt_m)=1}}
\frac{1}{t_dt_m}\sum_{\substack{t,k  \\ tk=t_dt_m }}t^{2g} \mu \pare{ \frac{t_dt_m}{t} }
\mu(t_d)\Co_{\wh{\tau} }\pare{\frac{q^{\frac{n^2}{2}}\H_{\tau'}(q) }{q-1} }^{2(g-1)}\nonumber\\
& = & \sum_{\substack{ \tau,\wh \tau ,t_d,t_m \\ |\tau|=n \\ \wh{\tau} = \tau/(t_d,t_m) \\ (\dhat,t_dt_m)=1}}
\frac{\mu(t_d)}{t_dt_m}\Co_{\wh{\tau} }\pare{\frac{q^{\frac{n^2}{2}}\H_{\tau'}(q) }{q-1} }^{2(g-1)}
\sum_{\substack{t,k  \\ tk=t_dt_m }}t^{2g} \mu \pare{ \frac{t_dt_m}{t} }\nonumber
\end{eqnarray}
which, thanks to \eqref{f:torus} gives:
\begin{equation}
\NSMF_n(q) =  \frac{1}{(q-1)^{2(g-1)}}
\sum_{\substack{ \tau,\wh \tau , t_d,t_m \\ |\tau|=n \\ \wh{\tau} = \tau/(t_d,t_m) \\ (\dhat,t_dt_m)=1}}
\mu(t_d)\frac{O^{2g}(t_dt_m )}{t_dt_m}\Co_{\wh{\tau} }\pare{q^{\frac{n^2}{2}}\H_{\tau'}(q) }^{2(g-1)}\label{f:C1toNSM}
\end{equation}
finishing thus the proof of formula \eqref{f:mainf}.
\endproof

\section{Corollaries} 
\index{Corollaries%
@\emph{Corollaries}}%


From formula \eqref{f:C1toNSM} of last section we can compute $\NSMF_n(1)$, which by Remark \ref{r:Echar}, 
computes the Euler Characteristic of $\MM(\Sl_n).$

\begin{rem}
As in 3.7 of \cite{VilHaus}, the only type $\tau$ whose $q^{\frac{n^2}{2}}\H_{\tau'}(q) $ has a simple zero 
at $q=1$ is the one of degree $d=n,$ and hence multiplicity $m_{n,(1)}=1$ and $0$ for other pairs 
$(d,\l)\ne (n,(1)).$ Let us write $\tau_n$ for the unique type of degree $n$. 
All other types have at least a double zero at $q=1$. 
\end{rem}

\begin{rem}
The value of 
$${\frac{q^{\frac{n^2}{2}}\H_{\tau_n'}(q) }{(q-1)}}$$ 
at $q=1$ is $n$. 
\end{rem}

\begin{cor}\label{c:Echar}
For $g\geq2$, the Euler Characteristic of $\MM(\Sl_n)$ is $\mu(n)n^{4g-3}$.
\end{cor}
\proof 
By Remark \eqref{r:Echar} we need to evaluate $\NSMF_n(1)$, which thanks to last two remarks gets simplified to
\begin{eqnarray}
\NSMF_n(1) & = & \sum_{\substack{ \tau,\wh \tau,t_d,t_m \\ \wh{\tau} = \tau/(t_d,t_m) \\ (\dhat,t_dt_m)=1}}
\mu(t_d)\frac{O^{2g}(t_dt_m )}{t_dt_m}\Co_{\wh{\tau} }\pare{\frac{q^{\frac{n^2}{2}}\H_{\tau'}(q) }{(q-1)}}^{2(g-1)}\Bigg|_{q=1}\nonumber\\
& = & \sum_{\substack{ \tau_n,\wh \tau_n,t_d,t_m \\ \wh{\tau}_n = \tau_n/(t_d,t_m) \\ (\dhat,t_dt_m)=1}}
\mu(t_d)\frac{O^{2g}(t_dt_m )}{t_dt_m}\Co_{\wh{\tau}_n }n^{2(g-1)}.\label{e:EChar1} 
\end{eqnarray}
Since the only multiplicity of $\tau_n$ is $1$, the only possible value for $t_m$ is $1$.

We shall now re-index the summation in \eqref{e:EChar1} by $t=t_d$ a divisor of $n$. And since the only possibles quotient types are those
with only one nonzero multiplicity $m_{\dhat,(1)}=1$, we can describe them by their degree, which is $n/t$ for these cases.

Therefore, \eqref{e:EChar1} becomes
\begin{eqnarray}
\NSMF_n(1) & = & \sum_{\substack{ t|n \\ (n/t,t)=1}}
\mu(t)\frac{O^{2g}(t)}{t}\Co_{\wh{\tau}_n }n^{2(g-1)}\nonumber\\
& = & \sum_{\substack{ t|n \\ (n/t,t)=1}}
\mu(t)\frac{O^{2g}(t)}{t}\mu \pare{\frac{n}{t}}\frac{t}{n} n^{2(g-1)}\nonumber\\
& = & \sum_{\substack{ t|n \\ (n/t,t)=1}}
\mu(t)\mu \pare{\frac{n}{t}}  {O^{2g}(t)} n^{2g-3}\label{e:EChar2} 
\end{eqnarray}
and since $\mu(a)\mu(b)=\mu(ab)$ for relatively prime numbers $a$ and $b$, we have
\begin{equation}\label{e:Echar3}
\NSMF_n(1) = \sum_{\substack{ t|n \\ (n/t,t)=1}} \mu(n)  {O^{2g}(t)} n^{2g-3}.
\end{equation}
which, in turn, reduces to
\begin{equation}\label{e:Echar4}
\NSMF_n(1) = \mu(n) n^{2g-3} \sum_{ t|n } {O^{2g}(t)},
\end{equation}
since for square-free $n$, every divisor $t$ is relatively prime to $n/t$, 
and when $n$ is not square-free, all the summation is killed by $\mu(n)=0$.

Finally, the factor $\displaystyle{ \sum_{ t|n } {O^{2g}(t)}}$ is nothing but the total number of elements in the $2g$ torus 
$(\C^\times)^{2g}$ killed by $n$, namely $n^{2g}$ (see Remark \ref{r:torus}), concluding
\begin{equation}\label{e:Echar5}
\NSMF_n(1) = \mu(n) n^{2g-3} n^{2g}= \mu(n)n ^{4g-3}.
\end{equation}
\endproof

Another interesting fact about $E(q; {\MM}(\Sl_n))$ deduced from formula \eqref{e:NSM5} is the following 
\begin{cor}\label{c:palin}
The $E$-polynomials $E_n(q):=E(q; {\MM}(\Sl_n))$ of the Character Varieties $\MM(\Sl_n)$ are palindromic. Thanks to
Remark \ref{r:degepoly} 
this is equivalent to
\begin{equation}\label{e:Epalin}
E_n(q)= E_n(q^{-1})q^{2(n^2-1)(g-1)}.
\end{equation}
\end{cor}

\proof
This follows from the identity \ref{e:palin} we recall here:
$$\H_{\l}(q)=(-1)^{|\l|}\H_{\l'}(q^{-1}).$$ 
Since, the exponent $2(g-1)$ kills every power of $-1$ and we end up with
\begin{eqnarray}
E_n(q^{-1})q^{2(n^2-1)(g-1)} & = & 
\sum_{\tau,t}q^{2(n^2-1)(g-1)}\pare{q^{-\frac{n^2}{2}}\frac{\H_{\tau'}(q^{-1})}{q^{-1} - 1}}^{2(g-1)}t^{2g-1}C_\tau^t \nonumber \\
 & = & \sum_{\tau,t}\pare{q^{n^2-1}q^{-\frac{n^2}{2}}\frac{\H_{\tau'}(q^{-1})}{q^{-1} - 1} }^{2(g-1)}t^{2g-1}C_\tau^t \nonumber \\
 & = & \sum_{\tau,t}\pare{q^{\frac{n^2}{2}}\frac{\H_{\tau'}(q)}{1 - q} }^{2(g-1)}t^{2g-1}C_\tau^t \nonumber \\
 & = & \sum_{\tau,t}\pare{q^{\frac{n^2}{2}}\frac{\H_{\tau'}(q)}{q-1} }^{2(g-1)}t^{2g-1}C_{\tau'}^t \nonumber \\
 & = & E_n(q).\label{e:palinfinale} 
\end{eqnarray}
since $C_\tau^t$ and $C_{\tau'}^t$ agree by Remark \ref{r:tautau}.
\endproof

\chapter{The Kitchen}\label{ch:Kitchen}
\index{The Kitchen%
@\emph{The Kitchen}}%

In the next three sections we are going prove the remaining claims from Chapter \ref{ch:epoly}.

Lemma \ref{l:pm1} proves the assertion of Remark \ref{r:deltafactor} about \eqref{e:ZtoS}.

Lemmas \ref{l:whS} and \ref{l:mixnewdeg} take care of the statement made on Remark \ref{r:claimShat} and the one
made right after equation \eqref{e:Sform}, respectively.
Lemma \ref{l:mixdeg} proves the assertion made in Remark \ref{r:Zhatf}. 

Finally, in Section \ref{s:values} we prove Lemma \ref{l:keylemma}, finishing with the calculation of
the $\wh Z$'s and therefore with the proof of Theorem \ref{t:main}.

\section{The Sign}
\index{The Sign%
@\emph{The Sign}}%

Here we see how many values the exponential factor in \eqref{e:ZtoS} can take, and find sufficient conditions to get rid of it.

\begin{lem}\label{l:pm1}
The factor $\delta(\zeta_n)^{ds{\that\choose 2}{\sum|\Lambda(\gamma)|}}$ from \eqref{e:deltafactor} is $\pm 1.$ 
If $q$ satisfies \eqref{oddity}, then it is precisely $1$.
\end{lem}
\proof
The exponent ${ds{\that\choose 2}{\sum|\L(\gamma)|}}={ds{\that\choose 2}{\sum \l_i}}$ could be written as
\begin{eqnarray*}
{ds{\that\choose 2}{\sum \l_i}} &=& {\frac{1}{2}sd{\that (\that-1)}{\sum \l_i}} \\
&=&\frac{1}{2} s t\dhat (\that-1){\sum \l_i} \\
&=&\frac{1}{2} (q-1)\dhat (\that-1){\sum \l_i}
\end{eqnarray*}
which is one half of an integer divisible by $n,$ since $n|q-1$.

So the square of $\delta(\zeta_n)^{ds{\that\choose 2}{\sum|\Lambda(\gamma)|}}$ is
$\delta(\zeta_n)^{(q-1)\dhat (\that-1){\sum \l_i}}$ 
which is $1$, leading us to the first assertion of this Lemma, 
namely $$\delta(\zeta_n)^{ds{\that\choose 2}{\sum|\Lambda(\gamma)|}}=\pm 1.$$

For the second affirmation, we assume \eqref{oddity} and analyzing the exponent
$$\frac{1}{2} (q-1)\dhat (\that-1){\sum \l_i}$$ 
we see it is already divisible by $n.$ 

If $2n$ divides $q-1,$ then $n$ divides $\quottt{q-1}{2}$ and hence the whole exponent results divisible by $n.$

And when $n$ is odd, the exponent is an integer number whose double is divisible 
by $n$, implying that $n$ divides the exponent, already.
\endproof
\begin{rem}
We could also point out that since $\that | t$ and $t|n$, $\that -1$ must be even, so there is the extra factor of $2$ we were looking for. 
\end{rem}

\section{The Vanishing}
\index{The Vanishing%
@\emph{The Vanishing}}%

Let us now state and prove the vanishing lemmas. The proofs given here 
are rather elementary and follow from the same 
trick. 

\begin{lem}\label{l:whS}
$\wh{S}(t,\tau)=0$ if the type $\tau$ is not pure.
\end{lem}
\proof
The idea here is to twist all the lowest 
degree characters multiplying them by a fixed factor of $\alpha\in\Gamma_1,$ and get on one hand the original sum, 
but on the other one the same sum multiplied by a nontrivial factor, concluding thus that the whole sum is $0.$

More precisely, given the type $\tau=\{m_{d,\l}\}$, $d_1$ the lowest degree appearing in the support of $\tau$ 
(i.e.: the minimum $d$ such that $m_{d,\l}\ne 0$ for a nonempty $\l\in\PP$), and character $\chi=\chi_\L\in\Irr(G)$ 
of type $\tau$ fixed by the action of $\delta^s\in\Gamma_1$
we consider a new character $\chi^\alpha\in\Irr(G)$ corresponding to the function $\L^\alpha\PP^\Ji_n$ 
that maps every character in $\Gamma$ to the same partition that $\L$ does, except for those characters 
$\gamma\in\Gamma_{d_1}$ that are going to be mapped to $\L(\gamma\alpha^{-1}).$ It is immediate that $\chi^\alpha$ 
is also of type $\tau,$ fixed by the action of $\delta^s\in\Gamma_1$ and that 
$(\chi^\alpha)^{\alpha^{-1}}=\chi,$ so $\chi \mapsto \chi^{\alpha}$ is actually a bijection and hence the twisted 
summation 
$$\sum_{\substack{\chi\in \Irr(G)\\ \delta^s\chi=\chi\\ \tau(\chi)=\tau}}\Delta_{\chi^\alpha}(\zeta_n)$$ 
is nothing but a permutation of the terms involved in \eqref{e:Shat}, the sum 
for $\wh{S}(t,\tau),$ and therefore leads to the same value.

But 
\begin{eqnarray} 
\Delta_{\chi^\alpha}(\zeta_n) & = & \prod_{\{\gamma\}, \gamma \in \Gamma}\gamma^{|\Lambda^\alpha(\gamma)|}(\zeta_n) \\
& = & \prod_{\{\gamma\}, \gamma \in \Gamma_{d_1}}\gamma^{|\Lambda^\alpha(\gamma)|}(\zeta_n)
\prod_{\{\gamma\}, \gamma \notin \Gamma_{d_1}}\gamma^{|\Lambda^\alpha(\gamma)|}(\zeta_n) \nonumber\\
& = & \prod_{\{\gamma\}, \gamma \in \Gamma_{d_1}}\gamma^{|\Lambda(\gamma\alpha^{-1})|}(\zeta_n)
\prod_{\{\gamma\}, \gamma \notin \Gamma_{d_1}}\gamma^{|\Lambda(\gamma)|}(\zeta_n)\nonumber\\
& = & \prod_{\{\gamma\}, \gamma \in \Gamma_{d_1}}(\gamma\alpha)^{|\Lambda(\gamma)|}(\zeta_n)
\prod_{\{\gamma\}, \gamma \notin \Gamma_{d_1}}\gamma^{|\Lambda(\gamma)|}(\zeta_n)\nonumber
\end{eqnarray}
which, factoring the $\alpha$'s out becomes
\begin{eqnarray} 
\Delta_{\chi^\alpha}(\zeta_n) & = & \alpha(N_1^{d_1}(\zeta_n))^{\sum '|\Lambda(\gamma)|}
\prod_{\{\gamma\}, \gamma \in \Gamma}\gamma^{|\Lambda(\gamma)|}(\zeta_n) \label{e:exp1}\\
& = & \Delta_{\chi}(\zeta_n)\alpha(\zeta_n)^{\sum '|\Lambda(\gamma)|d_1} \label{e:exp2}\\
& = & \Delta_{\chi}(\zeta_n)\alpha(\zeta_n)^{\sum '|\l|d_1m_{d_1,\l}} \label{e:exp3}
\end{eqnarray}
where the summations in the exponents of \eqref{e:exp1}, \eqref{e:exp2} and \eqref{e:exp3} correspond respectively to:
\begin{itemize}
\item The sum of $|\Lambda(\gamma)|$ for every orbit  $\{\gamma\}$ with $\gamma$ of degree $d_1$
\item The sum of $|\Lambda(\gamma)|d_1$ for every orbit  $\{\gamma\}$ with $\gamma$ of degree $d_1$
\item The sum of $|\l|d_1m_{d_1,\l}$ for every partition $\l \in \PP.$
\end{itemize}

Going back to the sum \eqref{e:Shat} defining $\wh{S}(t,\tau)$ we conclude:
\begin{eqnarray} 
 \wh{S}(t,\tau) & = & \sum_{\substack{\chi\in \Irr(G)\\ \delta^s\chi=\chi\\ \tau(\chi)=\tau}}\Delta_{\chi^\alpha}(\zeta_n) \nonumber\\ 
 & = & \sum_{\substack{\chi\in \Irr(G)\\ \delta^s\chi=\chi\\ \tau(\chi)=\tau}}
 \Delta_{\chi}(\zeta_n)\alpha(\zeta_n)^{\sum |\l|d_1m_{d_1,\l}}\nonumber\\  
 & = & \wh{S}(t,\tau)\alpha(\zeta_n)^{\sum |\l|d_1m_{d_1,\l}}.\label{e:Shattwist}
\end{eqnarray}
If the exponent $\sum |\l|d_1m_{d_1,\l}$ is not $n$, the twisting factor 
$\alpha(\zeta_n)^{\sum |\l|d_1m_{d_1,\l}}$ could be nontrivial, since $n$ is the order of $\zeta_n.$ 
For instance, taking $\alpha=\delta$ (or any other generator of $\Gamma_1$) does the job.

So $\wh{S}(t,\tau)=0$ unless the exponent $\displaystyle{\sum_{\l\in\PP} |\l|d_1m_{d_1,\l}}$ is $n$. 
But $n$ is already the whole sum $\displaystyle{\sum_{d,\l} |\l|dm_{d,\l}}$. Thus the only way to have a 
nonzero $\wh{S}(t,\tau)$ is for $\tau$ to be pure of degree $\d_1$ , that is to say, all 
the $\gamma\in\Gamma$ involved should have the same degree $d=d_1$.
\endproof

\begin{lem}\label{l:mixnewdeg}
$\displaystyle{\wh{S}(t,\tau)=\sum_{\dhat|d}S(t,\dhat,\tau)}$, that is to say, we could 
disregard all those terms $\Delta_{\chi}(\zeta_n)$ corresponding to characters $\chi$ with mixed newdegrees.
\end{lem}
\proof
A similar trick works here. Assuming $\tau$ is a type pure of degree $d$, 
we define the new-type $\wh \tau$ of a type $\tau$ character $\chi_\L=\delta^s\chi_\L$ and the map $\L\in\PP^\Ji_n$ 
associated to it, as the collection of new-multiplicities $\wh m _{\dhat, \l}$ given by the number of Frobenius 
orbits $\{[\gamma]\}$ for $[\gamma]\in \quott{\Gamma}{\bra{\delta^s}}$, of newdegree $\dhat$ that are mapped to $\l\in\PP$ 
by $\L$. Let us write $\wh \tau \subseteq \tau.$

Then 
\begin{equation}\label{e:accShat}
\wh{S}(t,\tau)=\sum_{\wh \tau \subseteq \tau}S(t,\wh \tau, \tau) 
\end{equation}
where $S(t,\wh \tau, \tau)$ is the sum of $\Delta_{\chi}(\zeta_n)$ for those $\chi$ fixed by the action of $\delta^s$, of type $\tau$ and new-type $\wh \tau$.

Our claim will follow once we prove $S(t,\wh \tau, \tau)$ vanishes unless $\wh \tau$ is pure of some newdegree $\dhat$, and noting that
$S(t,\dhat,\tau)$ is the sum of all such $S(t,\wh \tau, \tau)$.

For the sums $S(t,\wh \tau, \tau)$ we proceed in a similar way as before, twisting by a factor of $\alpha\in\Gamma_1$ 
all those $\gamma$ of lowest new-degree $\dhat_1.$ Since multiplication by $\alpha$ preserves both degrees and newdegrees
it is easy to see that this twisting preserves both type and new-type.

In this case, from \eqref{e:accShat} we can pull out the factor
\begin{equation}\label{e:pullfact1}
\alpha(\zeta_n)^{\sum |\l|t\dhat_1\wh{m}_{\dhat_1,\l}} 
\end{equation}
which could be made nontrivial for a suitable choice of $\alpha$ ($\alpha=\delta$ works again) unless the exponent $\sum |\l|t\dhat_1\wh{m}_{\dhat_1,\l}$ is $n$, and this can only happen when all the $\gamma$ in the supports of the $\L$'s have the same newdegree $\dhat_1$, in other words, $\wh \tau$ is pure of newdegree $\dhat_1$.
\endproof

The claim from Remark \ref{r:Zhatf} will follow from next lemma.
\begin{lem}\label{l:mixdeg}

Let $n$ be a positive integer. Assume $\dhat,$ $t$ and a list $\vec{\l}=(\l_1,\ldots,\l_m)$ of positive integers are given such that
$n=\dhat t \sum \l_i.$ Let $q$ be a fixed prime power congruent to $1$ modulo $n,$ and note $s=(q-1)/t.$ 
Let $D$ be a multiple of $\dhat$ dividing $\dhat t.$ Write $\that$ for $\dhat t /D.$
Let $\Gamma_D$ and $\Gamma$ be as in Definition \ref{d:Gam}.
Consider $\vec{d}$ a list of degrees $(d_1,\ldots, d_m)$ all of which are multiples of $\dhat$ and divide $D.$

Then the sum
\begin{equation}\label{e:Zdeformated}
Z(\vec{d},t,{\dhat},\vec{\l}):=\sum_{\substack{([\alpha_1],\ldots,[\alpha_m])\in \pare{\quott{\Gamma_D}{\bra{\delta^s}}}^m\\ \{[\alpha_i]\}\neq\{[\alpha_j]\},\: i\neq j \\ \deg(\alpha_i)=d_i \\ \deg([\alpha_i])=\dhat}}\prod_{i=1}^m\alpha_i^{\l_i}(\zeta_n^{\that})
\end{equation}
will vanish if not all the entries from $\vec{d}$ are the same (the evaluations in the product being computed as in $\Gamma_D$).
\end{lem}

\proof
Let $S$ be the subset of $\pare{\quott{\Gamma_D}{\bra{\delta^s}}}^m$ consisting on all the $m$-tuples of character classes
$([\alpha_1],\ldots,[\alpha_m])$ as in the summation, namely $[\alpha_i]\in\quott{\Gamma_D}{\bra{\delta^s}}$ has degree 
$\dhat$ and $\alpha_i$ has degree $d_i,$ no two of the entries $[\alpha_i]$ being $Frob_q$-conjugates.

Since every character in $\Gamma_1$ has degree $1,$ multiplication by $\beta\in\Gamma_1$ changes neither degree nor newdegree.
Then, we can define a free $\quott{\Gamma_1}{\bra{\delta^s}}$ action on $S$ as follows.
A class character $[\beta]\in\quott{\Gamma_1}{\bra{\delta^s}}$ acts on 
the $m$-tuple $([\alpha_1],\ldots,[\alpha_m])$ by coordinatewise multiplication on those entries $[\alpha_i]$ correspond 
to those $i$ with $d_i=\min\set{d_j}$ (it is easy to see that in the new $m$-tuple no two of the entries are $Frob_q$-conjugates).

We could regard $\pare{\quott{\Gamma_D}{\bra{\delta^s}}}^m$ as the dual of an $m$-fold product of copies of $\Ker({}^TN^D_1(\delta^s))$ 
\footnote{i.e.: the subgroup of $\F_{q^D}^\times$ whose dual is $\quott{\Gamma_D}{\bra{\delta^s}}$ by Remark \ref{r:subquotdual}.}
, and ${\quott{\Gamma_1}{\bra{\delta^s}}}$ as a the subgroup of $\pare{\quott{\Gamma_D}{\bra{\delta^s}}}^m$ of $m$-tuples
$([\beta_1],\ldots,[\beta_m])$ with trivial entries on those $i$ with $d_i>\min \set{d_j}$ and the same character class 
$[\beta]$ in all others.

The result will follow by direct application of Lemma \ref{l:chthy} from Appendix, once we prove that 
the element $(\zeta_n^{\that},\ldots,\zeta_n^{\that})\in\Ker({}^TN^D_1(\delta^s))^m$ at which we are evaluating all the 
summation does not belong to the intersection of kernels from the statement.

And this is indeed true since, regarding $[\beta]$ as an $m$-tuple, evaluation at $\zeta_n^{\that}$ becomes
$$
[\beta](\zeta_n^{\that})= \beta(\zeta_n^{\that})^{D\sum' \l_i } = \beta(\zeta_n)^{{\that}D\sum' \l_i } 
$$ 
where the sum $\sum' \l_i$ from the exponents is over those $\l_i$ with $d_i=\min\set{d_j}.$
Taking $\beta=\delta$ and keeping in mind that $\zeta_n$ is of order $n$ and that not all of the $d_i$ are the same we get a nontrivial value
for  $[\beta](\zeta_n^{\that}).$
\endproof

\section{The Values}\label{s:values}
\index{The Values%
@\emph{The Values}}%

This whole section is devoted to the proof of Lemma \ref{l:keylemma}, namely
that the sums $\wh{Z}(D,t,\dhat,\vec{\lambda})$ defined in \eqref{e:Zhatf} satisfy
$$
\wh{Z}(D,t,\dhat,\vec{\lambda})=\left\{
\begin{array}{cc}
\mu(\dhat)(-\dhat)^{m-1}(m-1)!\frac{q-1}{t} & \text{when } \gcd(D,t)=1  \\
0 & \text{otherwise}  \end{array}
\right .
$$
Let us remind the definition of $\wh{Z}(D,t,\dhat,\vec{\lambda})$
$$
\wh{Z}(D,t,\dhat,\vec{\lambda})=\sum_{\substack{([\alpha_1],\ldots,[\alpha_m])\in (\quott{\Gamma_D}{\langle\delta^s\rangle})^m\\ \{[\alpha_i]\}\neq\{[\alpha_j]\},\: i\neq j \\ \deg([\alpha_i])=\dhat}}\prod_{i=1}^m\alpha_i^{\lambda_i}(\zeta_n^{\that})
$$
where we are assuming
\begin{itemize}
\item $q$ is a fixed prime power satisfying \eqref{oddity}, 
\item $D$, $\dhat$, $t$, $\that$ and $n$ are positive integers, 
\item $\dhat$ is a divisor of $D$,
\item $\dhat t = D \that $ and
\item $\vec{\l}=(\l_1,\ldots,\l_m)$ an $m$-tuple of positive integers satisfying 
\begin{equation} \label{e:cond}
n=D \;\that\; \sum_{i=1}^m \l_i, 
\end{equation}
\end{itemize}
and the conventions are that $s$ stands for $(q-1)/t,$
$\delta$ is some fixed generator of $\Gamma_1= \dual {\F}_q^\times,$ and
$\zeta_n\in\F_q^\times$ an element of order $n.$ 

\begin{rem}
The sum is taken over all possible $m$-tuples $([\alpha_1],\ldots,[\alpha_m])$ of classes of characters 
$[\alpha_i]\in \quott{\Gamma_D}{\bra{\delta^s}}$ of newdegree exactly $\dhat$ with the additional 
condition that no two of the $[\alpha_i]$ can be conjugate by the action of Frobenius, in particular,
they must be different. In other words, the orbits $\{[\alpha_i]\}$ and $\{[\alpha_j]\}$ must be 
disjoint, and since they are orbits, it is the same to say that they are different.
\end{rem}
\begin{rem}
The $\alpha_i$ are evaluated as in $\Gamma_D.$
\end{rem}

The outline of the proof goes as follows:
\begin{itemize}
\item[\bf first:] realize $\wh Z$ as one value of a particular function $f$ on certain poset $\Pi^\rho$ 
of set-partitions fixed by an action.

\item[\bf second:] define the accumulated sums $\wh f$ of $f.$ 

\item[\bf third:] compute $\wh f$.

\item[\bf fourth:] solve for $f$ using M\"obius Inversion Formula for $\wh f$ on $\Pi^\rho$.
\end{itemize}

\subsection{First step: The poset $\Pi^\rho.$}
\index{First step: The poset $\Pi^\rho.$%
@\emph{First step: The poset $\Pi^\rho.$}}%

Let us consider the poset $\Pi^\rho=\Pi^\rho_{m\dhat}$ (as in example \ref{ex:setp} on \ref{s:MIV}) of 
set-partitions of $\{1,2,\ldots, m\dhat\}$ fixed by the permutation $\rho$ consisting of $m$ 
disjoint $\dhat$-cycles
\begin{equation}\label{e;action}
\rho = (1,\ldots ,\dhat)(\dhat+1,\ldots ,2\dhat)\ldots ((m-1)\dhat+1,\ldots, m\dhat)\in \S_{m\dhat}.
\end{equation}

To every such set-partition $\nu$ assign $f(\nu)$ defined as the following sum: 

The range of the summation will be the set $R(\nu)$ of any possible ordered $m\dhat$-tuple 
$\vec \alpha = ([\alpha_{j}])_{i=1}^{m\dhat}$ of classes of characters in $\quott{\Gamma_D}{\bra{\delta^s}}$ satisfying both
\begin{itemize}\label{d:Rnu}
 \item[1] $[\alpha_{\rho(i)}]=[\alpha_{i}]^q $, so $\rho$ acts as Frobenius.
 \item[2] $[\alpha_i]=[\alpha_j]$ if and only if $i$ and $j$ belong to the same part of $\nu$, so the 
set-partition $\nu$ defines all the equality relations among the classes of characters in the $m\dhat$-tuple.
\end{itemize}
And the number we sum for $\vec \alpha \in R(\nu)$ is computed as
\begin{equation}
\prod_{ \substack{i,j \\ 1\leq i\leq m \\ j=(i-1)\dhat +1 } }\alpha_{j}^{\lambda_i}(\zeta_n^{\that})
\end{equation}
where each evaluation is made in $\Gamma_D$.


\begin{rem}
A couple of things should be checked for this definition to make sense.

The evaluation $\alpha(\zeta_n^{\that})$ in $\Gamma_D$ is well defined, that is to say, it does not depend on the representative $\alpha$ 
chosen from the class $[\alpha]\in\quott{\Gamma_D}{\bra{\delta^s}}$.

To see this we must check that both $\alpha(\zeta_n^{\that})$ and $(\delta^s\alpha)(\zeta_n^{\that})$ 
have the same value when the evaluation is done in $\Gamma_D.$ According to Remark \ref{r:subquotdual} this is the same 
that $\zeta_n^{\that}\in\Ker({}^TN^D_1(\delta^s)).$ 

This is indeed true since:
$$(\delta^s\alpha)(\zeta_n^{\that}) = \delta^s(N_1^D(\zeta_n^{\that}))\alpha(\zeta_n^{\that})
= \delta(\zeta_n^{sD\that})\alpha(\zeta_n^{\that})=\alpha(\zeta_n^{\that})$$
because $\zeta_n\in\F_q$, hence $N_D(\zeta_n)=\zeta_n^D$, and since
$$sD\that=st\dhat=(q-1)\dhat$$ 
which is divisible by $n,$ then
$\delta(\zeta_n^{sD\that})=\delta(1)=1.$ 
\end{rem}

With this setting, $\wh{Z}(D,t,\dhat,\vec{\lambda})$ becomes:
$$f(\wh 0)=\sum_{\vec \alpha \in R(\wh 0)}
\prod_{\substack{i,j \\ 1\leq i\leq m \\ j=(i-1)\dhat +1}}\alpha_{j}^{\lambda_i}(\zeta_n^{\that})$$ for $\wh 0$ 
the finest set partition (i.e. that with $m\dhat$ parts of one element each).

\subsection{Second step: The accumulated sums $\wh f$}
\index{Second step: The accumulated sums $\wh f$%
@\emph{Second step: The accumulated sums $\wh f$}}%

The good thing about this setting is that it is fairly easy to compute the accumulated sums for $f$ defined as
\begin{equation}\label{e:accf}
\wh f (\nu):=\sum_{\varrho\geq\nu} f(\varrho) 
\end{equation}
where the restrictions given by the set partition $\varrho$ are not necessarily strict, in other words, 
there could be repetitions among the classes of characters corresponding to different parts of $\varrho$.

More precisely 
\begin{equation}\label{e:fhat}
\wh f (\nu)=\sum_{\vec \alpha \in \wh R (\nu)}
\prod_{\substack{i,j \\ 1\leq i\leq m \\ j=(i-1)\dhat +1}}
\alpha_{j}^{\l_i}(\zeta_n^{\that}) 
\end{equation}
where 
\begin{equation}\label{e:Rhat}
\wh R(\nu) := \bigcup_{\varrho\geq\nu} R(\varrho)
\end{equation}
is defined analogously to $R(\nu)$ with the second condition in \eqref{d:Rnu} relaxed to:
\begin{itemize}
 \item[2'] $[\alpha_i]=[\alpha_j]$ if $i$ and $j$ belong to the same part of $\nu$.
\end{itemize}


\subsection{Third Step: Calculating $\wh f$}
\index{Third Step: Calculating $\wh f$%
@\emph{Third Step: Calculating $\wh f$}}%

Our next task is to prove that most of the $\wh f (\nu)$ are zero except, possibly, for the top one, 
the one corresponding to the partition $\wh 1$ having only one part. 

In \ref{sss:Kill} we prove this by pulling out a factor in the sum \eqref{e:fhat}. Such factor resulting a full character
sum of a cyclic group. To prove that the element at which the characters are being evaluated is 
nontrivial we count the size of certain torsion subgroup by computing $\gcd$'s. 

The calculation of $\wh f (\wh 1)$ is finished in \ref{sss:lo} when the size of the aforementioned cyclic group is finally computed.

\subsubsection{The Killing Factor}\label{sss:Kill}
\index{The Killing Factor%
@\emph{The Killing Factor}}%


To compute $\wh f (\nu)$ for a particular $\nu\in\Pi^\rho$, let us consider $\nu_1$ the part of $\nu$ containing the index $1.$

Let us define $\nu_s$ by 
$$\nu_s:=\bigcup_{i}\rho^i(\nu_1)$$
the union of $\rho$-conjugates from $\nu_1.$

Since the conditions on the classes of characters $[\alpha_i]$ makes them to coincide when their indices 
belong to the same part of $\nu,$ but there is no condition making them necessarily different when their indices 
belong to different parts, we can pull out a factor of 
\begin{equation}\label{e:facfhat}
\wh f_1 (\nu) :=\sum_{ \wh R _1(\nu) } \prod_{\substack{i,j \\ 1\leq i\leq m \\ j=(i-1)\dhat +1\\j \in \nu_s}} \alpha_{j}^{\l_i}(\zeta_n^{\that})
\end{equation}
from $\wh f (\nu)$, where the range of summation $\wh R _1(\nu)$ consists in all possible tuples of classes 
of characters $[\alpha_i],$ with the set of indices restricted to the union of $\rho$-orbits of elements in $\nu_1$
with the same conditions:
\begin{itemize}
 \item $[\alpha_{\rho(i)}]=[\alpha_{i}]^q $,
so $\rho$ acts as Frobenius.
 \item $[\alpha_i]=[\alpha_j]$ if $i$ and $j$ belong to the same part of $\nu$. 
\end{itemize}

Our next task is to prove that such factor $\wh f_1 (\nu)$ will vanish (making the whole $\wh f (\nu)$ vanish as well) 
unless, possibly, when $\nu_1$ is the whole partition $\nu.$ Equivalently, for $\wh f (\nu)$ to be nonzero we must have
$\nu=\wh 1$ the partition with only one part.

\begin{rem}
We must stress that since $\zeta_n\in\F_q,$ the product in the summation \eqref{e:facfhat} defining $\wh f_1 (\nu)$
could be regarded as a product of repeated classes of characters. This is because
$$\alpha^q(\zeta_n^{\that})=\alpha(\zeta_n^{q\that})=\alpha(\zeta_n^{\that})$$
and by definition of $\nu_1$ the $\alpha_i$ in the product are $Frob_q$-conjugates. 

Hence
\begin{equation}\label{e:facfsimp}
\prod_{\substack{i,j \\ 1\leq i\leq m \\ j=(i-1)\dhat +1\\j \in \nu_s}} 
\alpha_{j}^{\l_i}(\zeta_n^{\that}) = \alpha_{1}^{\sum '\l_i}(\zeta_n^{\that})
\end{equation}
where the exponent $\sum '\l_i$ is the sum of those $\l_i$ with $1\leq i\leq m$ such that $(i-1)\dhat +1 \in \nu_s.$ 
\end{rem}

\begin{rem}\label{r:degcharclass}
The presence of a subindex $k\in \nu_1$ greater than $1$ but smaller than $\dhat +1$ makes the class $[\alpha_1]$ 
to have a Frobenius orbit of a smaller size than $\dhat,$ that is to say $[\alpha_1]^{q^{k-1}}=[\alpha_k]=[\alpha_1]$. 
\end{rem}

Thus, the factor $\wh f_1 (\nu)$ end up being:
\begin{equation}\label{e:sumsimp}
\wh f_1 (\nu) = \sum_{ \wh R _1(\nu) } \alpha_{1}(\zeta_n^{\that\sum '\l_i})
\end{equation}
where the $[\alpha_1]$ ranges over 
the subgroup of $\quott{ \Gamma_D}{ \bra{\delta^s} }$ of classes of characters 
fixed by $Frob_q^{\dhat_1},$ 
where $\dhat_1$ is the number of parts $\nu_j$ of $\nu$ containing a $\rho$-conjugate of 
$1.$ \footnote{in other words, $\dhat_1$ is the degree of $\nu_1$ under the induced $\rho$-action.}

\begin{rem}
If $\dhat_1=\dhat$ then all the $i\in\nu_1$ are congruent to $1$ modulo $\dhat$, and the character class 
$[\alpha_1]$ ranges over all 
$[\alpha]\in\quott{\Gamma_D}{ \bra{ \delta^s }} $ 
satisfying $[\alpha]^{ q^{\dhat} }= [ \alpha ]$. 
\end{rem}
\begin{rem}\label{r:dhat11}
If $\dhat_1=1$, we must have $\{1,2,\ldots, \dhat \} \subseteq \nu_1$ and 
$[\alpha_1]$ moves on the subgroup of newdegree $1$ character classes of $\quott{\Gamma_D}{\bra{\delta^s}}$, namely those with
$[\alpha]^q=[\alpha]$. 
\end{rem}

\begin{rem}\label{r:weirdgroup}
In \eqref{e:sumsimp} we are computing a character sum (as in Remark \ref{r:partialcharsum}) over 
\begin{equation}
\pare{\quott{\Gamma_D}{\bra{\delta^s}} }_{\dhat_1}:=
\set{[\alpha]\in \pare{\quott{\Gamma_D}{\bra{\delta^s}} }\:/\: [\alpha]=[\alpha]^{q^{\dhat_1}}}
\end{equation}
a subgroup of a quotient of the group of characters $\Gamma_D,$
and we are evaluating them at $\zeta_n^{\that\sum '\l_i},$ with the summation in the exponent as explained in
\eqref{e:facfsimp}.
\end{rem}

\begin{rem}\label{r:deltahere}
Since the character $\delta\in\Gamma_1$ has degree $1,$ it must also have newdegree $1$ (by \ref{r:dhatdividesd}), and
hence it belongs to $\pare{\quott{\Gamma_D}{\bra{\delta^s}} }_{\dhat_1}$. 
\end{rem}

\begin{rem}\label{r:weirdgroup2}
The fact that we are considering the quotient $\quott{\Gamma_D}{\bra{\delta^s}}$ amounts of taking the character group
of the subgroup of $\F_{q^D}^\times$ killed by $\delta^s$ (by Remark \ref{r:subquotdual}), namely 
\begin{equation}
\quott{\Gamma_D}{\bra{\delta^s}}=\Hom\pare{\Ker\pare{ {}^TN_1^D(\delta^s)},\C^\times},
\end{equation}
and we have already checked that $\zeta_n^{\that}$ belongs to such kernel.
\end{rem}

\begin{rem}\label{r:weirdgroup3}
The fact that the sum is taken over a subgroup of a character group corresponds to take a full 
character sum on a quotient of such group (also by Remark \ref{r:subquotdual}),
which amounts to identify those elements in $\Ker\pare{{}^TN_1^D(\delta^s)}$ 
whose quotient belongs to the 
intersection $K$ of the kernels of the $[\alpha]$'s involved in the sum \eqref{e:sumsimp}.
In other words
\begin{equation}\label{e:Kdef}
K := \bigcap_{[\alpha]\in\pare{\quott{\Gamma_D}{\bra{\delta^s}} }_{\dhat_1}}\Ker([\alpha]) 
\:\subseteq \: \Ker\pare{{}^TN_1^D(\delta^s)} \:\subseteq\: \F_{q^D}^\times
\end{equation}
and  
\begin{equation}\label{e:predual}
\pare{\quott{\Gamma_D}{\bra{\delta^s}} }_{\dhat_1} =  
\Hom\pare{{\quott{\Ker\pare{{}^TN_1^D(\delta^s)}}{K}},\C^\times}.
\end{equation}
\end{rem}

By Remarks \ref{r:weirdgroup}, \ref{r:weirdgroup2} and \ref{r:weirdgroup3}  we conclude the following 

\begin{lem}\label{l:inKornot}
The character sum \eqref{e:sumsimp}
will be zero if $\zeta_n^{\that\sum '\l_i}$ does not belong to the intersection of kernels $K,$ 
And if it does, the result will be the order of the group
$\pare{\quott{\Gamma_D}{\bra{\delta^s}} }_{\dhat_1}.$
\end{lem}

Since ${\Gamma_D}$ is cyclic, then so is $\quott{\Gamma_D}{\bra{\delta^s}}$ and its subgroup 
$\pare{\quott{\Gamma_D}{\bra{\delta^s}} }_{\dhat_1}$. 

The order of $\quott{\Gamma_D}{\bra{\delta^s}}$ is 
$$\frac{|\Gamma_D|}{|\bra{\delta^s}|}= \frac{q^D-1}{t}$$
and in a cyclic group of such order, the subgroup of elements fixed by $Frob_q^{\dhat_1}$ 
(which by definition is that of those elements killed by $q^{\dhat_1}-1$) agrees with the one of elements killed by 
\begin{equation}\label{e:order}
\gcd \pare{\frac{q^D-1}{t}, q^{\dhat_1}-1} 
\end{equation}
and thus this must also be its order.

\begin{lem}\label{l:ordweirdgroup}
The order of $\pare{\quott{\Gamma_D}{\bra{\delta^s}} }_{\dhat_1}$  is
$$\frac{q^{\dhat_1}-1}{t} \gcd \pare{ \frac{D}{\dhat_1} , t }.$$
\end{lem}

\proof
Keeping in mind that $t|n$ and $n|q-1$ 
we could simplify \eqref{e:order} to
\begin{eqnarray}
\gcd \pare{\frac{q^D-1}{t}, q^{\dhat_1}-1}  & = & 
\frac{q^{\dhat_1}-1}{t} \gcd \pare{ \underbrace{\pare{q^{\dhat_1}}^{\frac{D}{\dhat_1} -1}+\ldots+q^{\dhat_1}+1}_{\text{ $D/\dhat_1$ terms }} , t }
\nonumber \\
&=& 
\frac{q^{\dhat_1}-1}{t} \gcd \pare{ \underbrace{1 +\ldots+1+1}_{\text{ $D/\dhat_1$ terms }} , t } \nonumber\\
&=& 
\frac{q^{\dhat_1}-1}{t} \gcd \pare{ \frac{D}{\dhat_1} , t }.\label{e:order2}
\end{eqnarray}
The Lemma is thus proved.
\endproof

\begin{rem}\label{r:inKornot}
Plugging this result in Lemma \ref{l:inKornot}, we conclude that the value of the character sum \eqref{e:sumsimp}
\begin{equation}
\wh f_1 (\nu) = \sum_{ \wh R _1(\nu) } \alpha_{1}(\zeta_n^{\that\sum '\l_i}) 
\end{equation}
will be either 
$$\frac{q^{\dhat_1}-1}{t} \gcd \pare{ \frac{D}{\dhat_1} , t }$$
or $0$ depending on whether
$\zeta_n^{\that\sum '\l_i}$ belongs or not to the intersection $K$ from \eqref{e:Kdef}.
\end{rem}

We are now in conditions to prove the following
\begin{lem}\label{l:nonzero} 
All the accumulated sums $\wh f (\nu)$ vanish except, possibly, for the top one $\wh f (\wh 1).$
\end{lem}

\proof
Let us assume $\wh{f}(\nu)\ne 0.$

Since we factored out $\wh{f}_1(\nu)$ from the sum \eqref{e:fhat}, the vanishing of former implies that of the latter.

From Lemma \ref{l:inKornot} we know $\wh{f}_1(\nu)$ will be zero unless
$\zeta_n^{\that\sum '\l_i}\in K.$

$K$ was defined as an intersection of kernels in \eqref{e:Kdef}, and
one of the character classes involved in such intersection is $[\delta]$ 
(Remark \ref{r:deltahere}). 

Let us see that $\zeta_n^{\that\sum '\l_i}$ does not belong to $\Ker([\delta])$ unless
the sum $\sum'\l_i$ in the exponent (as defined in \eqref{e:facfsimp}) is the whole sum 
$\sum \l_i = \dfrac{n}{D\that}.$

Evaluation of $\delta\in\Gamma_D$ at $\zeta_n^{\that\sum '\l_i}$ gives:
$$
\delta(N_1^D(\zeta_n^{\that\sum '\l_i}))=\delta(\zeta_n^{D\that\sum '\l_i})
$$
(evaluated as in $\Gamma_1$) and for this to be $1$, we must have $n| D\that\sum '\l_i$ because $\zeta_n$ 
is of order $n$ and $\delta$ generates $\dual \F_q^\times$.

Since $n= D\that\sum \l_i$ (see assumption \eqref{e:cond}) we conclude 
$$ \sideset{}{'}\sum \l_i= \sum \l_i = \frac{n}{D\that},$$
which means all the exponents $\l_i$ are involved in the sum $\sum '.$ 

According to \eqref{e:facfsimp} and the definition of $\wh{R}_1(\nu)$ (right after \eqref{e:facfhat}), in order to 
have $\zeta_n^{\that\sum '\l_i}$ killed by $\delta$ (viewed in $\Gamma_D$) 
we must have 
\begin{equation}\label{e:wholesumpartition}
\nu = \bigcup_{i}\rho^i(\nu_1)
\end{equation}
and thus
\begin{equation}\label{e:wholesumindex}
\wh R (\nu) = \wh R _1(\nu)
\end{equation}
therefore
$\displaystyle{\sum_{ \wh R _1(\nu) } \alpha_{1}^{\sum '\l_i}(\zeta_n^{\that}) }$ is not just a factor
of \eqref{e:fhat}, but the whole sum 
$$ \wh f (\nu)=\sum_{\vec \alpha \in \wh R (\nu)}
\prod_{\substack{i,j \\ 1\leq i\leq m \\ j=(i-1)\dhat +1\\j \in \nu}}
\alpha_{j}^{\l_i}(\zeta_n^{\that}).$$

So far we know that $\wh{f}_1(\nu)$ will vanish unless \eqref{e:wholesumpartition} holds, and in that case
$\wh{f}_1(\nu) = \wh{f}(\nu).$ 
Let us see that it also vanish unless $\nu=\nu_1,$ that is to say, $\nu=\wh 1.$

Taking a closer look at 
$$\zeta_n^{\that\sum '\l_i}=\zeta_n^{\that\sum \l_i} = \zeta_n^{\frac{n}{D}}=\zeta_D$$
we see that the element at which we are evaluating the character classes in \eqref{e:sumsimp}
is an order $D$ element.

According to Lemma \ref{l:inKornot}, in order to compute $\wh f_1 (\nu)$ 
we only need to check whether $\zeta_D$ belongs or not to the intersection $K$ 
from \eqref{e:Kdef}, and since $K$ is cyclic, this is equivalent as checking whether $D$ divides or not 
the order of $K$ (Remark \ref{r:group2}).

Analyzing formula \eqref{e:predual} from Remark \ref{r:weirdgroup3} and Remark \ref{r:weirdgroup2}, 
we see that the order of $K$ is the quotient of the orders of $\quott{\Gamma_D}{\bra{\delta^s}}$ and
$\pare{\quott{\Gamma_D}{\bra{\delta^s}} }_{\dhat_1},$ which by 
Lemma \ref{l:ordweirdgroup} is the quotient
\begin{equation}\label{e:ordK}
\frac{\frac{q^D-1}{t}}{\frac{q^{\dhat_1}-1}{t} \gcd \pare{ \frac{D}{\dhat_1} , t }}=
\frac{q^D-1}{(q^{\dhat_1}-1)\gcd \pare{ \frac{D}{\dhat_1} , t }}
\end{equation}
and an order $D$ element like $\zeta_D$ will belong to it if and only if $D$ divides its order (Remark \ref{r:group}). 

For $D$ to divide $\displaystyle{\frac{q^D-1}{(q^{\dhat_1}-1)\gcd \pare{ \frac{D}{\dhat_1} , t }} }$ it is necessary that it divides
$$\frac{q^D-1}{q^{\dhat_1}-1}=\pare{q^{\dhat_1}}^{\frac{D}{\dhat_1} -1}+\ldots+q^{\dhat_1}+1 \equiv \frac{D}{\dhat_1} \mod(q-1)$$
but $D|n$ and $n|q-1$. So if $D$ divides $\frac{D}{\dhat_1}$ that means $\dhat_1=1.$

By Remark \ref{r:dhat11}, and \eqref{e:wholesumpartition} we conclude that the only way that $D$ divides
the order of $K$ (and hence $\wh{f}(\nu)$ could be nonzero) is if 
$$\nu = \bigcup_{i}\rho^i(\nu_1)= \nu_1$$
meaning $\nu$ is the top set-partition $\wh 1$, the one with only one part.
\endproof

\subsubsection{The Leftovers}\label{sss:lo}
\index{The Leftovers%
@\emph{The Leftovers}}%

So far we know that the only potentially nonzero accumulated sum $\wh f (\nu)$ is the top one $\wh f (\wh 1) = f(\wh 1),$ therefore
$f(\wh 0)=\mu(\wh 0,\wh 1) \wh f(\wh 1)$ by M\"obius Inversion Formula.

According to Remark \ref{r:weirdgroup}, Lemma \ref{l:inKornot}, 
and \eqref{e:ordK} from the proof of Lemma \ref{l:nonzero}, 
$\wh f(\wh 1)$ could still be either \footnote{the order of $\pare{\quott{\Gamma_D}{\bra{\delta^s}}}_1$.}
$\displaystyle{\frac{q-1}{t} \gcd \pare{ {D} , t }}$
or zero, depending on whether $D$ divides or not
$\displaystyle{ \frac{q^D-1}{(q-1)\gcd \pare{ {D} , t }} }.$

Let us further analyze the quotient $\dfrac{q^D-1}{q-1}$. 

By the Newton binomial expansion 
\begin{eqnarray*}
q^D-1 &=& (q-1 + 1)^D-1 \\
&=& -1 + 1 + (q-1)D + (q-1)^2{D \choose 2} + \sum_{k\geq 3}(q-1)^k{D \choose k}\\
&=& (q-1)D + (q-1)^2{D \choose 2} + O((q-1)^3) 
\end{eqnarray*}
where $O((q-1)^k)$ stands for an integer multiple of $(q-1)^k$.

Then, to get a nonzero $\wh{f}(\wh 1),$ $D$ must divide
\begin{eqnarray*}
\frac{q^D-1}{(q-1)\gcd \pare{ {D} , t }} &=&  
\frac{(q-1)D + (q-1)^2{D \choose 2} + O((q-1)^3) }{(q-1)\gcd \pare{ {D} , t }}  \\
&=& \frac{D + (q-1){D \choose 2} + O((q-1)^2) }{\gcd \pare{ {D} , t }}
\end{eqnarray*}
or, what is the same, ${\gcd \pare{ {D} , t }}$ must divide
\begin{eqnarray*}
 \frac{D + (q-1){D \choose 2} + O((q-1)^2) }{D}  &=& 1 + \frac{(q-1)(D-1)}{2} + \frac{(q-1)}{D}O(q-1).
\end{eqnarray*}
Keeping in mind that ${\gcd \pare{ {D} , t }} | D$ and $D|q-1$ we can disregard the last term in the right hand side 
and end up with the condition
${\gcd \pare{ {D} , t }}$ divides $1 + \dfrac{(q-1)(D-1)}{2}$.

Now, for $n$ odd, $D-1$ results even and $\dfrac{(q-1)(D-1)}{2}$ is a multiple of $D$, hence of ${\gcd \pare{ {D} , t }}$ and we get
${\gcd \pare{ {D} , t }}=1.$

For even $n$ but $2n|q-1$, the term $\dfrac{(q-1)}{2}(D-1)$ is again multiple of $n$ and hence of $D$, therefore
${\gcd \pare{ {D} , t }}$ is again $1$.

Summarizing, if $q$ satisfies \ref{oddity}, $\wh{f}(\wh 1)$ will be zero unless
${\gcd \pare{ {D} , t }}=1,$ in which case it will be the order of
the order of $\pare{\quott{\Gamma_D}{\bra{\delta^s}}}_1.$

\subsection{Fourth step: The Inversion}
\index{Fourth step: The Inversion%
@\emph{Fourth step: The Inversion}}%

The assertion will now follow from an inversion formula on the poset $\Pi^\rho$, 
and the result of Hanlon on the M\"obius function of this poset 
(\cite{Han} Theorem 4). 

Under the current conditions,
 the only nonzero $\wh f (\nu)$ is $\wh f (\wh 1)=f(\wh 1)$ and hence
$$f(\wh 0)=\mu(\wh 0,\wh 1) \wh f (\wh 1)=
\left\{
\begin{array}{cc}
\mu(\wh 0,\wh 1) \frac{q-1}{t} \gcd \pare{ {D} , t }& \text{ if } \gcd(D,t)=1  \\
0 & \text{otherwise}  \end{array}
\right.
$$
But 
\begin{equation}\label{e:Hanlon}
\mu(\wh 0,\wh 1)=\mu(\dhat)(-\dhat)^{m-1}(m-1)! 
\end{equation}
and the factor $\gcd \pare{ {D} , t }$ is trivially removed in the first case, so we finally proved 
Lemma \ref{l:keylemma}:
$$
Z(D,t,\dhat,\vec{\lambda})=
f(\wh 0)=
\left\{
\begin{array}{cc}
\mu(\dhat)(-\dhat)^{m-1}(m-1)!\frac{q-1}{t} & \text{ if } \gcd(D,t)=1  \\
0 & \text{otherwise}  \end{array}
\right.
$$


%
%

\appendix
\index{Appendix@\emph{Appendix}}%

\chapter{Two Lemmas from Group Theory}\label{ap:twolem}
\index{Appendix!Two Lemmas from Group Theory@\emph{Two Lemmas from Group Theory}}%

In this appendix we will state and prove two lemmas we used during the computation of the number of 
points of $\MM(\Sl_n(\F_q))$ that are group-theoretical in nature.

\begin{lem}\label{l:grthy}
Let $G$ and $H$ be two finite groups acting on a set $X$. Let us make the further assumption that both actions commute, that is to say
$g\cdot (h\cdot x)=h\cdot (g\cdot x)$ for all $g\in G$, $h\in H$ and $x\in X$. 
There is a natural $G$-action on the set of orbits $X/H$ given by $g\cdot Hx;= H(g\cdot x)$ 
\footnote{well defined since both actions commute}. 
Let us note $[x]$ the $H$-orbit of $x\in X$. Then for any $x\in X$, the size of the $G$-orbit of $[x]$ divides that of the $G$-orbit of $x$.
\end{lem}
\proof
The $G$-stabilizer of $x$ is a subgroup of the $G$-stabilizer of $[x]$, hence the index of the latter divides that of the former. 
The lemma follows since the indexes coincide with the sizes of the corresponding $G$-orbits.
\endproof

\begin{lem}\label{l:chthy}
Let $G$ be an abelian group, $\dual G =\Hom(G, \C^\times)$ its dual, $g$ and element of $G$, and $\dual Q \subseteq \dual G$ a subgroup.
Pointwise multiplication defines an action $\dual Q \acts \dual G,$ and consider $S$ a subset of $\dual G$ stable by this $\dual Q$-action, and 
the character sum:
\begin{equation}\label{e:charsumS}
S(g) :=\sum_{\phi\in S} \phi(g).
\end{equation}
If $g$ does not belong to $\displaystyle K=\bigcap_{\psi\in \dual Q} \Ker(\psi)$ then $S(g)=0.$
\end{lem}
\proof
Since $\dual Q$ is a subgroup of $\dual G,$ it must be (by Remark \ref{r:subquotdual}) the dual of some 
quotient $Q=\quottt{G}{K}$ of $G$, with $K$ as in the statement (i.e.: the intersection of kernels of 
the characters in $\dual Q$).

The fact that $g\notin K$ amounts for the class $[g]\in Q$ being nontrivial. 

Then, the character sum 
\begin{equation}\label{e:charsumlem}
\dual Q (g):=\sum_{\psi\in \dual Q} \psi(g) 
\end{equation}
vanishes, since it is a full character sum of a nontrivial element.

The Lemma follows then from the observation that the condition on $S$ being $\dual Q$-stable makes $S$ the disjoint union of 
$\dual Q$-orbits, and each of them is a multiple of \eqref{e:charsumlem}, which is zero under our hypothesis.
\endproof




\bibliographystyle{plain}  
\bibliography{RefsBIS}     

\begin{thebibliography}{10}

\bibitem{DelLus}
P.~Deligne and G.~Lusztig.
\newblock Representations of reductive groups over finite fields.
\newblock {\em Ann. of Math. (2)}, 103(1):103--161, 1976.

\bibitem{Del:2}
Pierre Deligne.
\newblock Th{\'e}orie de {H}odge. {II}.
\newblock {\em Inst. Hautes {\'E}tudes Sci. Publ. Math.}, (40):5--57, 1971.

\bibitem{Del:3}
Pierre Deligne.
\newblock Th{\'e}orie de {H}odge. {III}.
\newblock {\em Inst. Hautes {\'E}tudes Sci. Publ. Math.}, (44):5--77, 1974.

\bibitem{Fran}
Fran\c{c}ois Digne and Jean Michel.
\newblock {\em Representations of finite groups of {L}ie type}, volume~21 of
  {\em London Mathematical Society Student Texts}.
\newblock Cambridge University Press, Cambridge, 1991.

\bibitem{FdQ}
Daniel~S. Freed and Frank Quinn.
\newblock Chern-{S}imons theory with finite gauge group.
\newblock {\em Comm. Math. Phys.}, 156(3):435--472, 1993.

\bibitem{FulHar}
William Fulton and Joe Harris.
\newblock {\em Representation theory}, volume 129 of {\em Graduate Texts in
  Mathematics}.
\newblock Springer-Verlag, New York, 1991.
\newblock A first course, Readings in Mathematics.

\bibitem{Green:old}
J.~A. Green.
\newblock The characters of the finite general linear groups.
\newblock {\em Trans. Amer. Math. Soc.}, 80:402--447, 1955.

\bibitem{Han}
Phil Hanlon.
\newblock The fixed-point partition lattices.
\newblock {\em Pacific J. Math.}, 96(2):319--341, 1981.

\bibitem{VilHaus}
Tam{\'a}s Hausel and Fernando Rodriguez-Villegas.
\newblock Mixed {H}odge polynomials of character varieties.
\newblock {\em Invent. Math.}, 174(3):555--624, 2008.
\newblock With an appendix by Nicholas M. Katz.

\bibitem{Kac}
V.~G. Kac.
\newblock Infinite root systems, representations of graphs and invariant
  theory.
\newblock {\em Invent. Math.}, 56(1):57--92, 1980.

\bibitem{Green:new}
M.~T. Karkar and J.~A. Green.
\newblock A theorem on the restriction of group characters, and its application
  to the character theory of {${\rm SL}(n,\,q)$}.
\newblock {\em Math. Ann.}, 215:131--134, 1975.

\bibitem{Leh}
Gustav~Isaac Lehrer.
\newblock The characters of the finite special linear groups.
\newblock {\em J. Algebra}, 26:564--583, 1973.

\bibitem{MacD}
I.~G. Macdonald.
\newblock {\em Symmetric functions and {H}all polynomials}.
\newblock Oxford Mathematical Monographs. The Clarendon Press Oxford University
  Press, New York, second edition, 1995.
\newblock With contributions by A. Zelevinsky, Oxford Science Publications.

\bibitem{Serre}
Jean-Pierre Serre.
\newblock {\em Linear representations of finite groups}.
\newblock Springer-Verlag, New York, 1977.
\newblock Translated from the second French edition by Leonard L. Scott,
  Graduate Texts in Mathematics, Vol. 42.

\bibitem{TopGal}
Jean-Pierre Serre.
\newblock {\em Topics in {G}alois theory}, volume~1 of {\em Research Notes in
  Mathematics}.
\newblock Jones and Bartlett Publishers, Boston, MA, 1992.
\newblock Lecture notes prepared by Henri Damon [Henri Darmon], With a foreword
  by Darmon and the author.

\bibitem{Sesh}
C.~S. Seshadri.
\newblock Geometric reductivity over arbitrary base.
\newblock {\em Advances in Math.}, 26(3):225--274, 1977.

\bibitem{Stan}
Richard~P. Stanley.
\newblock {\em Enumerative combinatorics. {V}ol. 1}, volume~49 of {\em
  Cambridge Studies in Advanced Mathematics}.
\newblock Cambridge University Press, Cambridge, 1997.
\newblock With a foreword by Gian-Carlo Rota, Corrected reprint of the 1986
  original.

\bibitem{VinRy}
Nathaniel Thiem and C.~Ryan Vinroot.
\newblock On the characteristic map of finite unitary groups.
\newblock {\em Adv. Math.}, 210(2):707--732, 2007.

\end{thebibliography}
\index{Bibliography@\emph{Bibliography}}

\printindex



\end{document}